\documentclass[12pt]{article}

\usepackage{amstext    }
\usepackage{amsthm    }
\usepackage{a4}
\usepackage[mathscr]{eucal}
\usepackage{mathrsfs}

\usepackage{amsmath}
\usepackage{amssymb}
\usepackage{amscd}

\newtheorem{theorem}{Theorem}[section]
\newtheorem{definition}[theorem]{Definition}
\newtheorem{proposition}[theorem]{Proposition}
\newtheorem{corollary}[theorem]{Corollary}
\newtheorem{lemma}[theorem]{Lemma}
\newtheorem{remark}[theorem]{Remark}

\newtheorem{example}[theorem]{Example}

\newtheorem{conjecture}[theorem]{Conjecture}

\newcommand{\cali}[1]{\mathscr{#1}}

\newcommand{\volume}{{\rm vol}}

\newcommand{\supp}{{\rm supp}}
\newcommand{\const}{{\rm const}}
\newcommand{\dist}{{\rm dist\ \!}}

\newcommand{\loc}{{loc}}
\newcommand{\ddc}{dd^c}

\newcommand{\DSH}{{\rm DSH}}

\newcommand{\reg}{{\rm reg}}
\newcommand{\sing}{{\rm sing}}
\newcommand{\codim}{{\rm codim\ \!}}

\newcommand{\Cc}{\cali{C}}

\newcommand{\Ec}{\cali{E}}
\newcommand{\Fc}{\cali{F}}
\newcommand{\Gc}{\cali{G}}
\newcommand{\Hc}{\cali{H}}

\newcommand{\Oc}{\cali{O}}

\newcommand{\C}{\mathbb{C}}
\newcommand{\N}{\mathbb{N}}
\newcommand{\Z}{\mathbb{Z}}
\newcommand{\R}{\mathbb{R}}
\renewcommand{\P}{\mathbb{P}}

\newcommand{\E}{{\cal E}}
\renewcommand{\H}{{\cal H}}

\title{Equidistribution towards the Green current for holomorphic maps}
\author{Tien-Cuong Dinh and Nessim Sibony}

\begin{document}
\maketitle

\begin{abstract}
Let $f$ be a non-invertible holomorphic endomorphism of a projective
space and $f^n$ its iterate of order $n$.
We prove that the
pull-back by $f^n$ of a generic (in the Zariski sense)
hypersurface, properly normalized, converge to the Green current
associated to $f$ when $n$ tends to infinity. We also give an
analogous result for the pull-back of positive closed
$(1,1)$-currents and a similar result for regular
polynomial automorphisms of $\C^k$.
\end{abstract}




\section{Introduction} \label{section_introduction}

Let $f$ be a holomorphic endomorphism of algebraic degree $d\geq
2$ on the projective space $\P^k$. Let $\omega$ denote the
Fubini-Study form on $\P^k$ normalized so that $\omega$ is
cohomologous to a hyperplane or equivalently
$\int_{\P^k}\omega^k=1$. It is well-known that the sequence of
smooth positive closed $(1,1)$-forms $d^{-n} (f^n)^*(\omega)$
converges weakly to a positive closed $(1,1)$-current $T$ of mass
1. Moreover, $T$ has locally continuous potentials and is
{\it totally invariant}, i.e.
 $f^*(T)=dT$. We call $T$ {\it the Green
current} of $f$. The complement of the support of $T$ is the Fatou
set, i.e. the sequence $(f^n)$ is locally equicontinuous there. We
refer 
the reader to the survey \cite{Sibony} for background.
Our main results in this paper are the
following theorems, where $[\cdot]$ denotes the current of
integration on a complex variety.

\begin{theorem} \label{main_theorem}
Let $f$ be a holomorphic endomorphism of algebraic degree $d\geq
2$ of $\P^k$ and $T$ the Green current associated to $f$. There is
a proper analytic subset $\Ec$ of $\P^k$ such that if $H$ is a
hypersurface of degree $s$ in $\P^k$ which does not contain any
irreducible component of $\Ec$ then $d^{-n}(f^n)^*[H]$ converge to $sT$ in the
sense of currents
when $n$ tends to infinity. Moreover, $\Ec$ is totally invariant,
i.e. $f^{-1}(\Ec)=f(\Ec)=\Ec$.
\end{theorem} 

The exceptional set $\Ec$ will be explicitely constructed in Sections
\ref{section_exceptional} and \ref{section_cv}. It is the union of
totally invariant proper analytic subsets of $\P^k$ which are minimal.
That is, they have no proper analytic subsets which are totally
invariant, see Example \ref{example1}. That example shows that $\Ec$
is not the maximal totally invariant analytic set.
The previous result is in fact a consequence of the following one, see
also Theorem \ref{th_final} for a uniform convergence result.

\begin{theorem} \label{main_theorem_bis}
Let $f$ be a holomorphic endomorphism of algebraic degree $d\geq
2$ of $\P^k$ and $T$ the Green current associated to $f$. There is
a proper analytic subset $\Ec$ of $\P^k$, totally invariant, such
that if $S$ is a positive closed $(1,1)$-current of mass $1$ in
$\P^k$ whose local potentials are not identically $-\infty$ on any
irreducible component of $\Ec$ then $d^{-n}(f^n)^*(S)\rightarrow T$ as
$n\rightarrow\infty$. 
\end{theorem}

The space $\H_d$ of holomorphic maps $f$ of a given algebraic
degree $d\geq 2$ is an irreducible quasi-projective manifold. We will
also deduce from our study the following result due to Forn\ae ss
and the second author \cite{FornaessSibony2}, see also
\cite{FornaessSibony1, Sibony}.

\begin{theorem} \label{theorem_generic_map}
There is a dense Zariski open subset $\H_d^*$ of $\H_d$ such that
if $f$ is a map in $\H_d^*$ then $d^{-n}(f^n)^*(S_n)\rightarrow T$
for every sequence $(S_n)$ of positive closed $(1,1)$-currents of
mass $1$ in $\P^k$.
\end{theorem}

The rough idea in order to prove our main results is as follows. Write
$S=\ddc u+T$. Then, the invariance of $T$ implies that $d^{-n}
(f^n)^*(S) = d^{-n}\ddc (u\circ f^n) +T$. We have to show, in
different situations, that
$d^{-n}u\circ f^n$ converge to 0 in $L^1$. This implies that $d^{-n}
(f^n)^*(S)\rightarrow T$. 
So, we have to study
 the asymptotic contraction ({\it {\`a} la} Lojasiewicz)
by $f^n$. The main estimate is obtained using
geometric estimates and convergence
results for plurisubharmonic functions, 
see Theorem \ref{th_contraction}. 
If $d^{-n} u\circ f^n$ do not converge to 0,
then using that the possible contraction is limited,
we construct a limit $v$
with strictly positive Lelong numbers.
We then construct other functions $w_{-n}$ such that the
current $\ddc w_{-n}+T$ has Lelong numbers $\geq \alpha_0>0$ and
$w_0=d^{-n}w_{-n}\circ f^n$. It follows from the last identity that
$w_0$ has positive Lelong numbers on an infinite union of analytic
sets of a suitable dimension. The volume growth of these sets implies 
that the current associated to $w_0$ has too large self-intersection. This 
contradicts bounds due to
Demailly and M{\'e}o  \cite{Demailly1, Meo2}.
(One should notice that the Demailly-M{\'e}o estimates depend on the
$L^2$ estimates for the $\overline\partial$-equation; they were recently
extended to the case of compact K{\"a}hler manifolds by Vigny \cite{Vigny}).
The previous argument has to be applied inductively on totally invariant
sets for $f$, which are a priori singular and on which we inductively
show the convergence to 0, starting with sets of dimension 0. So, 
we also have to develop the basics of the theory of weakly
plurisubharmonic functions on singular analytic sets which 
is probably of independent interest. The advantage of this class of
functions is that it has good
compactness properties.

One may conjecture that totally invariant analytic sets should be
unions of linear subspaces of $\P^k$. The case of dimension $k=2$ is
proved in \cite{CerveauLinsNeto, SSU}. These authors
complete the result in \cite{FornaessSibony3}. If this were
true for $k\geq 3$, our proof
would be technically simpler. It is anyway interesting to carry the
analysis without any assumption on the totally invariant sets since our
approach may be extended to the case of meromorphic maps on compact
K{\"a}hler manifolds. At the end of the paper, we will consider the case
of regular polynomial automorphisms of $\C^k$.

The problem of convergence was first considered by Brolin for
polynomials in dimension 1 and then by Lyubich,
Freire, Lopes and Ma{\~n}{\'e} for rational maps in $\P^1$ \cite{Brolin, Lyubich, FLM}. 
In dimension $k=2$, Forn\ae ss and the second author proved that $\Ec$ is
empty when the local multiplicity of $f$ at every point is $\leq
d-1$, see \cite{FornaessSibony2}. This implies Theorem
\ref{theorem_generic_map} 
in dimension 2 for $S_n=S$. The proof in
\cite{FornaessSibony2} can be extended to the general case, see also \cite{Sibony}.

The family of hyperplanes in $\P^k$ is parametrized by a projective
space of dimension $k$. It follows from Theorem \ref{main_theorem} that
for a hyperplane $H$, generic in the Zariski sense, we have
$d^{-n}(f^n)^*[H]\rightarrow T$. Russakovskii and Shiffman have
proved this result for $H$ out of a pluripolar set in the
space of parameters \cite{RussakovskiiShiffman}. Analogous results for
subvarieties in arbitrary  K{\"a}hler manifolds were proved by the
authors in \cite{DinhSibony4}.
Concerning Theorems \ref{main_theorem} and \ref{main_theorem_bis},
our conditions are not optimal. Indeed, it might happen that the
potentials of $S$ are identically $-\infty$ on some components of $\Ec$
and still $d^{-n} (f^n)^*(S)\rightarrow T$.

In the case of dimension $k=2$, our results (except several uniform
convergences, e.g. Theorem \ref{th_final}) can be deduced from results
by Favre and Jonsson.
These authors say  that
their condition is necessary and sufficient in order to have the
previous convergence, see  \cite{FavreJonsson}, and they give needed
tools for the proof in   
\cite[p.310]{FavreJonsson2}. In which case,
 if the Lelong number of $S$ vanishes at generic points
on each irreducible component of an exceptional set then
$d^{-n}(f^n)^*(S)\rightarrow T$. The problem is still open in higher
dimension. 
When the Lelong number of $S$ is 0 at every point of
$\P^k$, the convergence $d^{-n}(f^n)^*(S)\rightarrow T$ was
obtained by Guedj  \cite{Guedj}, see also Corollary \ref{cor_lelong_zero}.
In these works, the problem of convergence is reduced to the study of sizes of
images of balls under iterates of $f$. This approach was first used in
\cite{FornaessSibony1, FornaessSibony2}.

Recall that the self-intersection 
$T^p:=T\wedge\ldots \wedge T$, $p$ times, defines a positive closed
$(p,p)$-current which is totally
invariant, i.e. 
$f^*(T^p)=d^pT^p$, see \cite{DinhSibony6} for the pull-back
operator on currents. It is natural to consider the analogous
equidistribution problem towards $T^p$.

\begin{conjecture}
Let $f$ be a holomorphic endomorphism of $\P^k$ of algebraic degree
$d\geq 2$ and $T$ its Green current. 
Then 
$d^{-pn}(f^n)^*[H]$ converge to $s T^p$ for every
analytic subset $H$ of $\P^k$ of pure codimension $p$ and of degree $s$
which is generic.
Here, $H$ is generic if 
either $H\cap E=\varnothing$ or $\codim H\cap E=p+\codim E$ 
for any irreducible component $E$ of every totally invariant analytic
subset of $\P^k$.
\end{conjecture}
We will see later that there are only finitely many analytic sets
which are totally invariant. 
Theorem \ref{main_theorem} proves the conjecture for $p=1$.
Indeed, in that case, it is equivalent to check the condition for $E$ minimal.
For $p=k$, the measure $\mu:=T^k$ is the unique
invariant measure of maximal entropy, see \cite{FornaessSibony2,
BriendDuval, Sibony}. 
In this case, the conjecture was proved by the authors in \cite{DinhSibony1}. 
Weaker results in this direction were obtained in \cite{FornaessSibony2} and
\cite{BriendDuval}. We will give some details in Theorem \ref{th_conjecture_k}.
For $2\leq p\leq k-1$, the authors have proved in
 \cite{DinhSibony7} that for $f$ in a Zariski dense open set
 $\H'_d\subset\H_d$, there is no proper analytic subset of $\P^k$
 which is totally invariant and that the conjecture holds. 
Indeed, a version of Theorem \ref{theorem_generic_map}
is proved.


\section{Plurisubharmonic  functions} \label{section_wpsh}

We refer the reader to \cite{Hormander, Demailly2, DinhSibony4} for the
basic properties of plurisubharmonic (psh for short) and quasi-psh functions
on smooth manifolds. 
In order to study the Levi problem for analytic spaces $X$, the psh functions
which are considered, are the restrictions of psh functions on an open set of $\C^k$ for a
local embedding of $X$. 
Let $u:X\rightarrow \R\cup\{-\infty\}$ be an upper semi-continuous
function which is not identically equal to $-\infty$ on any
irreducible component of $X$. 
Forn\ae ss-Narasimhan proved that if $u$ is
subharmonic or equal to $-\infty$ on any holomorphic disc in $X$, then
$u$ is psh in the above sense \cite{FornaessNarasimhan}. However, this
class does not satisfy good
compactness properties which are useful in our analysis. 
Assume that $X$ is an analytic space of pure dimension $p$. Let 
$\reg(X)$ and $\sing(X)$ denote the regular and the
singular parts of $X$. We consider the following weaker notion of psh
functions which is modeled after the notion of weakly holomorphic functions.
The class has good compactness properties.

\begin{definition} \label{def_psh} \rm
A function $v:X\rightarrow \R\cup\{-\infty\}$ is {\it wpsh} if 
\begin{enumerate}
\item[(a)] $v$ is psh on $X\setminus\sing(X)$.
\item[(b)] For $a\in\sing(X)$, $v(a)=\limsup v(x)$ with $x\in\reg(X)$ and
      $x\rightarrow a$.
\end{enumerate} 
\end{definition}

Forn\ae ss-Narasimhan's theorem implies that psh functions are wpsh. 
Wpsh functions are psh when $X$ is smooth. One should notice that the
restriction of a wpsh function to an irreducible component of $X$ is not
necessarily wpsh. For example, consider $X=\{xy=0\}$ in the unit ball of $\C^2$, let $v=0$ on
$\{x=0\}\setminus (0,0)$ and $v=1$ on $\{y=0\}$, then $v$ is wpsh on
$X$ but its restriction to $\{x=0\}$ is not wpsh. Consider the
(strongly) psh function $v_n:=|x|^{1/n}$ on $X$. The sequence
$v_n$ converge to $v$ in $L^1(X)$. So, psh
functions on analytic sets do not have good compactness properties.

\begin{proposition} \label{prop_psh_extension}
Let $Z\subset X$ be an analytic subset of dimension $\leq p-1$ and
$v'$ a wpsh function on $X\setminus Z$. If $v'$ is locally bounded
from above near $Z$ then there is a unique wpsh function $v$ on $X$
equal to $v'$ outside $Z$. 
\end{proposition}
\proof 
The extension to a psh function on $\reg(X)$ is well-known. So, we can
assume that $Z\subset \sing(X)$. 
Condition (b) in Definition
\ref{def_psh} implies the uniqueness of the extension of $v'$. Define $v(a)=\limsup v(x)$ with
$x\not\in Z$ and $x\rightarrow a$. It is clear that $v=v'$ out of $Z$
and $v$ satisfies the conditions in Definition \ref{def_psh}.
\endproof

Now assume for simplicity that $X$ is an analytic subset of pure
dimension $p$ of an open set $U$ in $\C^k$. The general case can be
deduced from this one.
The following results give 
characterizations of wpsh functions.

\begin{proposition} Let $\pi:\widetilde X\rightarrow X\subset U$ be a
  desingularization of $X$. If $v$ is a wpsh function on $X$ then
  there is a psh function $\widetilde v$ on $\widetilde X$ such that 
$v(x)=\max_{\pi^{-1}(x)}\widetilde v$ for $x\in X$. Conversely, if
$\widetilde v$ is psh on $\widetilde X$ then $x\mapsto
\max_{\pi^{-1}(x)}\widetilde v$ defines a wpsh function on $X$.
\end{proposition}
\proof
Define $\widetilde v:=v\circ \pi$ outside the analytic set
$\pi^{-1}(\sing(X))$. 
This function is psh and is locally bounded above
near $\pi^{-1}(\sing(X))$. We can extend it to a psh function on $\widetilde X$
that we also denote by $\widetilde v$.
For $x\in
X$, $\pi^{-1}(x)$ is compact. The maximum principle implies
that $\widetilde v$ is constant on each irreducible component of
$\pi^{-1}(x)$. From the definition of wpsh function, 
we get $v(x)=\max_{\pi^{-1}(x)}\widetilde
v$. The second assertion in the proposition follows from  the definition of
wpsh functions.
\endproof

A theorem of Lelong says that 
the integration on $\reg(X)$ defines a  positive closed
$(k-p,k-p)$-current $[X]$ on $U$, see \cite{Lelong, Demailly2}. Let $z$
denote the coordinates in $\C^k$. 

\begin{proposition} \label{prop_psh_equi}  
\rm 
A function $v:X\rightarrow \R\cup\{-\infty\}$ is {\it wpsh} if and only
if the following properties are satisfied :
\begin{enumerate}
\item[(a)] $v$ is in $L_\loc^1(X)$, i.e. $\int_K |v|
(\ddc|z|^2)^p<+\infty$ for any compact set $K\subset X$.

\item[(b)] $v$ is {\it strongly} upper semi-continuous, i.e. for
any $a\in X$ and any full measure subset $X'\subset X$ we have
$v(a)=\limsup v(x)$ with $x\in X'$ and $x\rightarrow a$.

\item[(c)] $\ddc(v[X])$ is a positive current on $U$.
\end{enumerate}
\end{proposition}
\proof
We use the notations in Proposition 2.3.
The proposition is known for smooth manifolds, see
\cite{Demailly2}. Assume that $v$ is wpsh. 
 The function $\widetilde v$ defined above satisfies properties (a), (b) and (c) on $\widetilde
X$. It follows that $v$ satisfies (a) and (b) on $X$. 
Since $\ddc (v[X])=\pi_*(\ddc(\widetilde v[\widetilde X]))$, $\ddc
(v[X])$ is positive. 
Hence, $v$ satisfies (c).

Conversely, Properties (a)-(c) imply that $v$ is psh on
$\reg(X)$. Then, Property (b) implies that $v$ satisfies the
conditions of
Definition \ref{def_psh}.
\endproof

\begin{proposition} \label{prop_psh_compactness}
Let $(v_n)$ be a sequence of wpsh functions on $X$, locally uniformly
bounded from above. Then,  there is a
subsequence $(v_{n_i})$ satisfying one of the following properties:
\begin{enumerate}
\item[{\rm (a)}] There is an irreducible  
component $Y$ of $X$ such that $(v_{n_i})$ converges uniformly to $-\infty$
on $K\setminus\sing(X)$ for any compact set $K\subset Y$.
\item[{\rm (b)}] $(v_{n_i})$ converges in $L^q_\loc(X)$ to a wpsh function $v$
for every $1\leq q<+\infty$. 
\end{enumerate}
In the last case, $\limsup
v_{n_i}\leq v$ on $X$ with equality almost everywhere.
\end{proposition}
\proof Let $\pi:\widetilde X\rightarrow X\subset U$ be as above.  We extend
the functions $v_n\circ \pi$, which are psh on $\pi^{-1}(\reg(X))$
to psh functions $\widetilde v_n$ on $\widetilde X$. 
Recall that $v_n(x)=\max_{\pi^{-1}(x)}\widetilde
v_n$. Now, since the proposition holds for smooth manifolds, it is
enough to apply it to $(\widetilde v_n)$. If a psh function
$\widetilde v$ is a limit value of $(\widetilde v_n)$ in
$L^q_\loc(\widetilde X)$, the function
$v$, defined by $v(x):=\max_{\pi^{-1}(x)}\widetilde v$, satisfies the
property (b) in the proposition. If not, $\widetilde v_n$ converge to
$-\infty$ locally uniformly on some component of $\widetilde X$ and
the property (a) holds.
\endproof

The following result is the classical {\it Hartogs' lemma} when $X$ is
smooth \cite{Hormander}.

\begin{lemma} \label{lemma_hartogs}
Let $(v_n)$ be a sequence of wpsh functions on $X$.
Let  $u$ be a continuous function on
$X$ such that $\limsup v_n<u$. Then for every compact set $K\subset
X$,  $v_n< u$ on $K$ for $n$ large enough.
This holds in particular, if $(v_n)$ converges 
to a wpsh function $v$ in $L_\loc^1(X)$ and $v<u$.
\end{lemma}
\proof
Let $\pi$ and $\widetilde v_n$ be defined as above. These functions
$\widetilde v_n$ are psh on $\widetilde X$. Define 
$\widetilde u:=u\circ \pi$. It is clear that $\widetilde u$ is 
continuous and  that $\limsup \widetilde v_n \leq \limsup  v_n\circ \pi <
\widetilde u$.
We only have to apply the classical Hartogs' lemma in order
to obtain $\widetilde v_n<\widetilde u$ on $\pi^{-1}(K)$ for $n$ large
enough. This implies the result. The last assertion in the lemma
is a consequence of Proposition \ref{prop_psh_compactness}.  
\endproof

The following lemma will be useful.

\begin{lemma} \label{lemma_restriction_compact}
Let $\Gc$ be a family of psh functions on $U$ 
locally  uniformly bounded from above. Assume that for each
irreducible component
of $X$ there is an analytic subset $Z$ such that the restriction of $\Gc$ to
  $Z$ is bounded in $L_\loc^1(Z)$. Then, the restriction of $\Gc$ to $X$ is
  bounded in $L^1_\loc(X)$.  
\end{lemma}
\proof
We can assume that $X$ is irreducible. For $(v_n)\subset \Gc$, define the
psh functions $\widetilde v_n$ on $\widetilde X$ as above. It is clear
that $\widetilde v_n$ are locally 
uniformly bounded from above. 
Let $W\Subset U$ be an open
set which intersects $Z$.
The maximal value of $\widetilde v_n$ on $\pi^{-1}(Z\cap \overline W)$ is equal
to the maximal value of $v_n$ on  $Z\cap \overline W$.
It follows from the hypothesis that no
  subsequence of $(\widetilde v_n)$ converges uniformly on compact
  sets to
  $-\infty$. Proposition \ref{prop_psh_compactness} 
applied to $(\widetilde v_n)$, 
implies that this sequence is
  bounded in $L^1_\loc(\widetilde X)$. Applying again Proposition \ref{prop_psh_compactness} to
  $(v_n)$ gives the lemma.  
\endproof

Let $R$ be a positive closed $(1,1)$-current on $U$ with
continuous local potentials, i.e. locally $R=\ddc v$ with $v$ psh and continuous.
Let  $R'$ be a positive closed
$(k-p,k-p)$-current on $U$, $1\leq p\leq k-1$. Recall that we can define
their intersection  by $R\wedge R':=\ddc(vR')$ where $v$ is
a local potential of $R$ as above. This is a positive closed
$(k-p+1,k-p+1)$-current on $U$ which depends continuously on $R'$. 
The definition is independent of the choice
of $v$. By induction, if $R_1$, $\ldots$, $R_p$ are positive closed
$(1,1)$-currents with continuous local potentials, 
the intersection $\nu:=R_1\wedge\ldots\wedge
R_p\wedge[X]$ is a positive measure with support in $X$.
This product is symmetric with respect to $R_1$, $\ldots$, $R_p$. 

\begin{proposition} \label{prop_cln}
For every compact sets $K$ and $K'$ with $K\Subset K'\subset X$, there is a constant $c>0$ such
that if $u$ is wpsh on $X$ we have
$$\max_Ku \leq c \|u\|_{L^1(K')} \quad \mbox{and} \quad \int_K
|u|d\nu\leq c \|u\|_{L^1(K')}.$$
In particular, $\nu$ has no mass on analytic subsets of dimension
$\leq p-1$ of $X$.  
\end{proposition}
\proof
Choose a compact set $L$ such that $K\Subset L\Subset K'$ and a
neighbourhood $W$ of $\sing(X)$ small enough.
If $a$ is a point in $K\cap W$, then  we can find a 
Riemann surface in $X$ containing $a$ and having boundary in
$L\setminus W$. Indeed, it is enough to consider the intersection of
$X$ with a suitable
linear plane $P$ of dimension $k-p+1$ passing throught $a$. The
maximum principle applied to the lift of $u$ to $\widetilde X$
(defined above) implies
that $u(a)\leq \max_{L\setminus W} u$ and hence 
$\max_K u \leq \max_{L\setminus W} u$. Since $L\setminus W\subset
\reg(X)$, the submean inequality for psh functions on smooth manifolds
implies that $\max_{L\setminus W} u\leq c\|u\|_{L^1(K')}$ for
some constant $c>0$. Hence, $\max_Ku \leq c \|u\|_{L^1(K')}$.

We prove now the second inequality. Replacing $u$ by $u-c\|u\|_{L^1(K')}$
allows us to assume that $u\leq 0$ on $K$. 
Since the problem is local, we can assume that $R_i=\ddc v_i$ with $v_i$
continuous on $U$. Moreover, we can approximate $v_i$ by decreasing
sequences $(v_{i,n})$ of smooth psh functions. Define $R_{i,n}:=\ddc
v_{i,n}$. It is well-known that
$\nu_n:=R_{1,n}\wedge \ldots\wedge R_{p,n}\wedge [X]$ converge
to $\nu$ in the sense of measures. Using the same arguments as in
the Chern-Levine-Nirenberg inequalities \cite{ChernLevineNirenberg, Demailly2, Sibony}
yields 
$$\int_K ud\nu_n\geq
-c'\|v_{1,n}\|_{L^\infty(K')}\ldots\|v_{p,n}\|_{L^\infty(K')} \|u\|_{L^1(K')}$$
where $c'>0$ is independent of $n$. When $n\rightarrow\infty$, since
$\nu_n\rightarrow\nu$ and since $u$ is upper semi-continuous, we
obtain
$$\int_K ud\nu\geq
-c'\|v_1\|_{L^\infty(K')}\ldots\|v_p\|_{L^\infty(K')} \|u\|_{L^1(K')}.$$
This implies the second inequality in the proposition.

Let $Y$ be an analytic subset of $X$ of dimension $\leq p-1$. Then,
there is a psh function $u'$ on $U$ such that
$\{u'=-\infty\}=Y$. The last inequality applied to the restriction of
$u'$ to $X$, implies $\nu(Y)=0$. 
\endproof


\section{Modulo $T$ plurisubharmonic functions}

We are going to develop in this section the analogue in the compact
case of the local theory in Section \ref{section_wpsh}.
Consider a (compact) analytic subset $X$ of $\P^k$ of pure dimension
$p$. Recall that the Green current $T$ of $f$ has locally continuous
potentials. Observe that in what follows (except for Lemma \ref{lemma_t_psh_interate},
Corollary \ref{cor_psh_convergence_T} and Remark \ref{rk_chaotic}), $T$ could be an arbitrary
positive closed $(1,1)$-current of mass 1 with continuous potentials, and
 $\P^k$ could be replaced by any compact K{\"a}hler manifold.
We will use the following notion that allows us to
simplify the notations.

\begin{definition} \rm
A function $u:X\rightarrow \R\cup\{-\infty\}$ is {\it wpsh modulo $T$}
if locally it is the difference of a wpsh function on $X$ and a potential of
$T$. If $X$ is smooth, we say that $u$ is {\it psh modulo $T$}. 
\end{definition}

The following result is a consequence of Proposition \ref{prop_psh_equi}.

\begin{proposition} \label{prop_equi_T_psh}
A function $u:X\rightarrow \R\cup\{-\infty\}$ is wpsh modulo $T$ if
and only if the following properties are satisfied
\begin{enumerate}
\item[{\rm (a)}] $u$ is in $L^1(X)$, i.e. $\int_X |u| \omega^p<+\infty$.

\item[{\rm (b)}] $u$ is {\it strongly} upper semi-continuous.

\item[{\rm (c)}] $\ddc(u[X])\geq - T \wedge [X]$ on $\P^k$.
\end{enumerate}
\end{proposition}

Note that if $u$ is a modulo $T$ wpsh function, 
$\ddc(u[X])+T\wedge [X]$ is a positive closed current of bidegree 
$(k-p+1,k-p+1)$ supported on $X$. If $S$ is a positive closed
$(1,1)$-current on $\P^k$ of mass $1$, then it is cohomologous to $T$ and 
we can write $S=T+\ddc u$ where $u$ is a modulo $T$ psh function on $\P^k$.
The restriction of such
a function $u$ to $X$
is either wpsh modulo $T$ or equal to $-\infty$ on at least one
irreducible component of $X$.

\medskip

The following proposition is a consequence
of Proposition \ref{prop_psh_compactness}.

\begin{proposition} \label{prop_T_psh_compactness}
Let $(u_n)$ be a sequence of modulo $T$ wpsh functions on $X$, uniformly
bounded from above. Then  there
is a subsequence $(u_{n_i})$ satisfying one of the following properties:
\begin{enumerate}
\item[{\rm (1)}] There is an irreducible component $Y$ of $X$ such that $(u_{n_i})$ converges uniformly to $-\infty$
on $Y\setminus\sing(X)$.
\item[{\rm (2)}] $(u_{n_i})$ converges in $L^q(X)$ to a modulo $T$ wpsh function $u$
for every $1\leq q<+\infty$. 
\end{enumerate} 
In the last case, $\limsup
u_{n_i}\leq u$ on $X$ with equality almost everywhere.
\end{proposition}

The Hartogs' lemma \ref{lemma_hartogs} implies the following.

\begin{lemma} \label{lemma_hartogs_T}
Let $(u_n)$ be a sequence of modulo $T$ wpsh functions on $X$ converging in
$L^1(X)$ to a modulo $T$ wpsh function $u$. If $w$ is a continuous function on
$X$ such that $u<w$, then $u_n< w$ for $n$ large enough.
\end{lemma}

The following lemma is deduced from Lemma \ref{lemma_restriction_compact}.

\begin{lemma} \label{lemma_restriction_compact_T}
Let $\Gc$ be a family of modulo $T$ psh functions on $\P^k$ 
uniformly bounded from above. Assume that each irreducible component of $X$
  contains an analytic subset $Y$ such that the restriction of $\Gc$ to
  $Y$ is bounded in $L^1(Y)$. Then, the restriction of $\Gc$ to $X$ is
  bounded in $L^1(X)$.  
\end{lemma}

Define a positive measure supported on $X$ by $\mu_X:=T^p\wedge [X]$. 
By B{\'e}zout's theorem, the mass of $\mu_X$ is equal to the degree of
$X$. The same argument implies that $\mu_X$ has positive mass on any
irreducible component of $X$.
The following result is a consequence of Proposition \ref{prop_cln}.

\begin{proposition} \label{prop_cln_mu}
There is a constant $c>0$ such that if $u$ is a
  modulo $T$ wpsh function on $X$ then
$$\max_Xu\leq  c (1+\|u\|_{L^1(X)}) \quad \mbox{and}\quad \int |u| d\mu_X\leq c (1+\|u\|_{L^1(X)}).$$
In particular, $\mu_X$ has no mass on analytic subsets of dimension
$\leq p-1$ of $X$.
\end{proposition}

We also have the following useful Proposition and Lemma.

\begin{proposition} \label{prop_cln_mu_bis}
A family $\Gc$ of modulo $T$ wpsh functions on $X$ is bounded
  in $L^1(X)$ if and only if there is a constant $c>0$ such that 
$|\int ud\mu_Y|\leq c$ for $u\in\Gc$ and for any irreducible component $Y$ of $X$.
\end{proposition}
\proof
Proposition \ref{prop_cln_mu} implies that
$\mu_Y$ has no mass on
$\sing(X)$. 
If $\Gc$ is bounded in $L^1(X)$ then it is bounded in $L^1(Y)$.
We have seen that the restriction of $u\in\Gc$ to $Y$ is equal outside
$\sing(X)$ to a modulo $T$ wpsh function on $Y$.
By Proposition \ref{prop_cln_mu},   there is a constant $c>0$ such that 
$|\int ud\mu_Y|\leq c$ for $u\in\Gc$.

Conversely, assume that $|\int ud\mu_Y|\leq c$ for $u\in\Gc$ and for
any irreducible component $Y$ of $X$. Since
$\mu_Y$ has no mass on
$\sing(X)$, we can replace $X$ by $Y$ and
assume that $X$ is irreducible.  Define $m_u:=\max_X u$ and
$v:=u-m_u$. 
Since $\max_X v=0$, Proposition \ref{prop_T_psh_compactness} implies that the family of
such functions $v$ is bounded in $L^1(X)$, see also Definition \ref{def_psh}(b).
On the other hand, we have
$$|m_u|\ \!\|\mu_X\|=\Big|\int ud\mu_X-\int v d\mu_X\Big|\leq c+\Big|\int vd\mu_X\Big|.$$ 
This and Proposition \ref{prop_cln_mu}, applied to $v$, imply that $|m_u|$ is
bounded.  Since $u=m_u+v$, we obtain that $\Gc$ is
bounded in $L^1(X)$.
\endproof

\begin{lemma} \label{lemma_t_psh_interate}
Let $u$ be a modulo $T$ wpsh function on $X$. 
If $X$ is invariant by $f$, i.e. $f(X)=X$, then $d^{-1}u\circ f$ is 
equal out of $\sing(X)\cup f^{-1}(\sing(X))$ to a modulo $T$ wpsh
function $w$ on $X$.
Moreover, $w$ depends continuously on $u$. 
\end{lemma}
\proof
Consider a point $x\in X$ out of $\sing(X)\cup f^{-1}(\sing(X))$. Since $T$
is totally invariant, if $v$ is
a potential of $T$ in a neighbourhood $V$ of $f(x)$ then $d^{-1}v\circ f$
is a potential of $T$ in a neighbourhood $U$ of $x$. Since the function $u+v$ is
psh on $X\cap V$, $d^{-1}(u\circ f+v\circ f)$ is psh on $X \cap U$. 
Hence, $\ddc ((d^{-1}u\circ f)[X])\geq -T\wedge [X]$ out of $\sing(X)\cup
f^{-1}(\sing(X))$. On the other hand, since $u$ is bounded from above, $d^{-1}u\circ f$
is bounded from above. Proposition \ref{prop_psh_extension} implies
the existence of $w$. That $w$ depends continuously on $u$ follows
from Proposition \ref{prop_T_psh_compactness}.
\endproof

\begin{corollary} \label{cor_psh_convergence_T}
Assume that $X$ is invariant. Let $\Gc$ be a
family of modulo $T$ wpsh functions on $X$, bounded in $L^1(X)$. Then, the family
  of  modulo $T$ wpsh functions on $X$ which are equal almost
  everywhere to $d^{-n} u\circ
  f^n$ with $n\geq 0$ and $u\in\Gc$, is
  bounded in $L^1(X)$. Moreover, if a modulo $T$ wpsh
  function $u$ on $X$ is a limit value of $(d^{-n} u_n\circ f^n)$ in
  $L^1(X)$ with $u_n\in\Gc$,
 then $u\leq 0$ on
  $X$ and $u=0$ on $\supp(\mu_X)$. The sequence $(d^{-n}u_n\circ f^n)$
  converges to $0$ in $L^1(\mu_X)$.
\end{corollary}
\proof
Replacing $f$ by an iterate $f^n$ allows us to
assume that $f$ fixes all the irreducible components of $X$. So,
we can assume that $X$ is irreducible.
For the first assertion,
by Propositions  \ref{prop_cln_mu},
we can subtract from each $u$ a constant in order that $\max_X
u=0$. So, we can assume that $\Gc$ is the set of such functions
$u$. This is a bounded set in $L^1(X)$. All the functions $d^{-n} u\circ
  f^n$ are equal almost everywhere to functions in $\Gc$. The first
assertion follows.

For the second assertion,
by Lemma \ref{lemma_t_psh_interate},
$d^{-n}u_n \circ f^n$ is equal outside an analytic set to a modulo $T$
wpsh function $v_n$ on $X$.
Propositions  \ref{prop_T_psh_compactness} and \ref{prop_cln_mu}  imply
that $u_n\leq A$ and $\int |u_n|d\mu_X\leq A$ for some constant
$A>0$. It follows that
$v_n\leq d^{-n}A$, see also Proposition \ref{prop_equi_T_psh}(b), and then 
$\limsup  v_n\leq 0$.
Hence, $u\leq 0$. On the other hand,
since $X$ is invariant and $T$ is totally invariant, we have
$(f^n)_*(\mu_X)=\mu_X$ and
$$\Big|\int (d^{-n} u_n\circ f^n) d\mu_X\Big| = d^{-n} \Big|\int u_nd (f^n)_*(\mu_X)\Big| = d^{-n} \Big|\int
u_n d\mu_X\Big|\leq d^{-n}A.$$
Hence, $\int v_nd\mu_X\rightarrow 0$. 
By Propositions
\ref{prop_cln_mu_bis} and \ref{prop_cln_mu}, 
$(v_n)$ is bounded from above. This allows us to apply the last
assertion in Proposition \ref{prop_T_psh_compactness}. We deduce from
Fatou's lemma and the convergence $\int v_nd\mu_X\rightarrow 0$, that
 $\int ud\mu_X\geq 0$. This and the inequality  
$u\leq 0$ imply that $u=0$ $\mu_X$-almost
everywhere. By upper semi-continuity, $u=0$ on $\supp(\mu_X)$.
\endproof

\begin{remark} \label{rk_chaotic} \rm
Assume that $f$ is chaotic, i.e. the support of the Green measure
$\mu$ of $f$ is equal to $\P^k$. Then, the previous corollary gives us
a simple proof of the following property: for all positive closed
$(1,1)$-currents $S_n$ of mass $1$ on $\P^k$, we have
$d^{-n}(f^n)^*(S_n)\rightarrow T$. Indeed, we can write $S_n=T+\ddc
u_n$ with $u_n$ bounded in $L^1(X)$, and hence $d^{-n} u_n\circ f^n$ converge to 0.
\end{remark}


\section{Lelong numbers} \label{section_lelong}

In this section, we recall some properties of the Lelong numbers of
currents and of plurisubharmonic functions, see \cite{Demailly2} for
a systematic exposition.

Let $R$ be a positive closed $(p,p)$-current on an open set $U$ of
$\C^k$. Let $z$ denote the coordinates in $\C^k$ and $B_a(r)$ the
ball of center $a$ and of radius $r$. Then, $R\wedge
(\ddc\|z\|^2)^{k-p}$ is a positive measure on $U$. 
Define for $a\in U$
$$\nu(R,a,r):=\frac{{\|R\wedge (\ddc\|z\|^2)^{k-p}\|_{B_a(r)}}}{{\pi^{k-p}r^{2(k-p)}}}.$$ 
When $r$ decreases to $0$, $\nu(R,a,r)$ is decreasing and
the {\it
Lelong number} of $R$ at $a$ is the limit
$$\nu(R,a):=\lim_{r\rightarrow 0} \nu(R,a,r).$$
The property that $\nu(R,a,r)$ is decreasing implies the 
following property that we will use later: {\it if $R_n\rightarrow R$ and
$a_n\rightarrow a$, then $\limsup \nu(R_n,a_n)\leq \nu(R,a)$}.

The Lelong number $\nu(R,a)$ is also the mass of the measure $R\wedge
(\ddc\log\|z-a\|)^{k-p}$ at $a$. It 
does not depend on
the coordinates. So, we can define the Lelong number for currents on
any manifold. If $R$ is the current of integration on an analytic set
$V$, by Thie's theorem, $\nu(R,a)$ is equal to
the multiplicity of $V$ at $a$. Recall also a theorem of Siu which
says that for $c>0$ the level set $\{\nu(R,a)\geq c\}$ is an analytic
subset of dimension $\leq k-p$ of $U$.

Let $S$ be a current of bidegree $(1,1)$ and $v$ a potential of
$S$ on $U$. Define {\it the Lelong number} of $v$ at $a$ by
$\nu(v,a):=\nu(S,a)$. We also have
\begin{equation}
\nu(v,a)=\lim_{r\rightarrow 0} \frac{\sup_{B_a(r)}v(z)}{\log r}.
\label{eq_lelong_psh}
\end{equation}
The function $\log r\mapsto \sup_{B_a(r)} v$ is increasing and
convex with respect to $\log r$. It follows that
if $v$ is defined on $B_a(1)$ and is negative, the fraction in
(\ref{eq_lelong_psh}) is decreasing when $r$ decreases to $0$. 
So, if two psh functions differ by a locally bounded function,
they have the same Lelong number at every point. Moreover, 
identity (\ref{eq_lelong_psh}) allows  to define the Lelong
number for every function which locally differs from a psh function by
a bounded function.

Let $X$ be an analytic subset of pure dimension $p$ in $U$ and $u$ a
wpsh function on $X$. Then,
$S^X:=\ddc (u[X])$ is a positive closed $(k-p+1,k-p+1)$-current on
$U$. Define 
$$\nu_X(u,a):=\nu(S^X,a).$$ 
When $X$
is smooth at $a$, we can also define a
positive closed $(1,1)$-current on a neighbourhood of $a$ in $X$ by  $S_X:=\ddc u$.  
We have $\nu_X(u,a)=\nu(S_X,a)$ where the last
Lelong number is defined on $X$.

Consider a proper finite holomorphic map $h:U'\rightarrow U$ between an
open set $U'$ of $\C^k$ and $U$. Let $X'$ be an analytic subset of
pure dimension $p$ of $U'$ such that
$h(X')=X$, and 
$a'\in U'$ a point such that $h(a')=a$.  It
follows from Proposition \ref{prop_psh_extension} that $u\circ h$ is
equal almost everywhere to a wpsh function $u'$
on $X'$. The continuity of $u'$ with respect to $u$ is proved as in Lemma \ref{lemma_t_psh_interate}.

\begin{proposition} \label{prop_lelong_compare}
Let $\delta$ denote the local topological degree of $h$ at $a'$. Then
$$\delta^{-k}\nu_X(u,a)\leq \nu_{X'}(u',a')\leq \delta\nu_X(u,a).$$
\end{proposition}
\proof 
Recall that $X$ and $X'$ may be reducible and singular, but one can
work on each irreducible component separately.
We deduce from the identity $h(X')=X$ and from the definition of $\delta$
that near $a:$ 
$$ \ddc(u[X]) \leq h_*(\ddc(u'[X']))\leq \delta
\ddc(u[X]).$$
Hence,
\begin{equation} \label{eq_lelong_compare}
\nu(\ddc(u[X]),a)\leq \nu(h_*(\ddc(u'[X'])),a)\leq \delta
\nu(\ddc(u[X]),a).
\end{equation}
On the other hand, by Theorems 9.9 and 9.12 in \cite{Demailly2},
we have
\begin{equation} \label{eq_lelong_compare_bis}
\nu(\ddc(u'[X']),a')\leq \nu(h_*(\ddc(u'[X'])),a) \leq \delta^k
\nu(\ddc(u'[X']),a').
\end{equation}
The inequalities in the proposition
follow from (\ref{eq_lelong_compare}) and
(\ref{eq_lelong_compare_bis}). 
\endproof

Let $B_a^X(r)$ denote the connected component of $B_a(r)\cap X$ which
contains $a$. We call it {\it the ball of center $a$ and of radius
  $r$ in $X$}.

\begin{proposition} \label{prop_psh_log}
Let $\Gc$ be a family of wpsh functions on $X$ which is compact in
$L^1_\loc(X)$. Let $\delta > 0$ such that
$\nu_X(u,a)< \delta$ for $u\in \Gc$ and $a\in X$.
Then, for 
any compact set $K\subset X$, there exist constants $c>0$ and $A>0$
such that
$$\sup_{B^X_a(r)}u\geq c\delta\log r -A \quad \mbox{ for }
u\in\Gc,\ a\in K  \mbox{ and } 0<r<1.$$
Moreover, the constant $c$ is independent of $\Gc$ and of $\delta$.
\end{proposition}
\proof 
Reducing $U$ allows to assume that $\Gc$ is bounded in $L^1(X)$ and $\nu_X(u,a)\leq \delta-\epsilon$,
$\epsilon>0$, on $X$ for every $u\in\Gc$. Moreover, 
by Proposition \ref{prop_cln}, $\Gc$ is 
uniformly bounded from above. So, we can assume that $u\leq 0$ 
for every $u\in\Gc$. If $0<r_0<1$ is fixed and $r_0<r<1$, the fact that
$\Gc$ is bounded in $L^1(X)$ implies that $\sup_{B^X_a(r)}u\geq
-A$ for every $a\in K$ where $A>0$ is a constant. Hence, it is enough to consider $r$
small.

We first consider the case where $X$ is smooth. 
Since the problem is local
we can assume that $X$ is a ball in $\C^p$. Up to a dilation of
coordinates, we can assume that the distance between $K$ and $\partial
X$ is larger than 1. 
Define 
$$s(u,a,r):=\frac{\sup_{B_a(r)\cap X} u}{\log r} .$$ 
Hence, for $a\in K$ and for $0<r<1$,
$s(u,a,r)$ decreases to $\nu(u,a)$ when $r$ decreases to 0. 
For every $(a,u)\in K\times\Gc$, since $\nu(u,a)\leq\delta-\epsilon$, 
there is an $r>0$ such that $s(u,a,r')\leq\delta -\epsilon/2$ for
$r'\leq 2r$. It follows that if a psh function $v$ on $X$ is close enough to $u$
then $s(v,a,r)\leq \delta-\epsilon/4$, see Lemma \ref{lemma_hartogs}. We then deduce
from the definition of $s(v,a,r)$ that if $b$ is close enough to $a$
and if $r'':=r-|b-a|$ then
$$s(v,b,r)\leq \frac{\log r''}{\log r} s(v,a,r'')\leq  
\frac{\log r''}{\log  r} s(v,a,r)\leq \delta.$$
The fact that $s(v,b,r)$ is increasing implies that 
$s(v,b,r')\leq \delta$ for $r'\leq r$ and for $(b,v)$ in a
neighbourhood of $(a,u)$. Since $K\times\Gc$ is compact, if $r$ is
small enough, the inequality $s(u,a,r)\leq\delta$ holds for
every $(a,u)\in K\times \Gc$. This implies the proposition for $c=1$
in the case where $X$ is smooth.

Now consider the general case. Since the problem is local, we can
assume that $X$ is analytic in $U=D_1\times D_2$ where $D_1$ and $D_2$
are the unit balls in
$\C^p$ and $\C^{k-p}$ respectively. We can also assume that the
canonical projection $\pi:D_1\times D_2\rightarrow D_1$ is proper
on $X$. Hence, $\pi:X\rightarrow D_1$ defines a ramified covering. Let
$m$ denote the degree of this covering.
For $u\in\Gc$, define a function $u''$ on $D_1$ by
\begin{equation}
u''(x):=\sum_{z\in \pi^{-1}(x)\cap X} u(z). \label{eq_psh_push}
\end{equation}
Since $\ddc
u''=\pi_*(\ddc (u[X]))\geq 0$, $u''$ is equal almost everywhere to
a psh function $u'$. 
It is easy to check that the family $\Gc'$ of these functions
$u'$  is compact in $L^1_\loc(D_1)$. Fix a ball $D$ containing
$\pi(K)$ such that $\overline D\subset D_1$. We need the following
Lojasiewicz type inequality, see \cite[Proposition
4.11]{FornaessSibony2}, which implies that $z\mapsto \pi^{-1}(z)\cap X$ is
H{\"o}lder continuous of exponent $1/m$ with respect to the Hausdorff
metric. The lemma is however more precise and is of independent interest.

\begin{lemma} \label{lemma_lojasiewicz}
There is a constant $A>0$ such that for $z\in D$ and
  $x\in X$ with $\pi(x)\in D$, we have
$$\dist( \pi^{-1}(z)\cap X,x)\leq A\dist(z, \pi(x))^{1/m}.$$
Moreover, if $y$ and $z$ are in $D$ we can write 
$$\pi^{-1}(y)\cap X=\{y^{(1)},\ldots, y^{(m)}\}\qquad \mbox{and} \qquad
\pi^{-1}(z)\cap X=\{z^{(1)},\ldots, z^{(m)}\}$$ 
so that 
$$\dist(y^{(i)},z^{(i)})\leq A\dist(y,z)^{1/m} \qquad\mbox{for } 1\leq i\leq m.$$
\end{lemma}
\proof
We prove the first assertion.
Let $x_j$, $p+1\leq j\leq k$, denote the last $k-p$ coordinates of
$x$. Let $z^{(1)}$, $\ldots$, $z^{(m)}$ denote the points in
$\pi^{-1}(z)\cap X$ and $z^{(1)}_{p+1}$, $\ldots$, $z^{(1)}_k$, $\ldots$,
$z^{(m)}_{p+1}$, $\ldots$, $z^{(m)}_k$  their last
$k-p$ coordinates. Here, the points in $\pi^{-1}(z)\cap X$ are repeated
according to their multiplicities. 
For $w\in D_1$, define $w^{(i)}$ and $w^{(i)}_j$ in the same way.
We consider the Weierstrass polynomials on $t\in\C$
$$\prod_{i=1}^m (t-w_j^{(i)})=t^m+a_{j,m-1}(w)
t^{m-1}+\cdots+a_{j,0}(w)=P_j(t,w).$$
The coefficients of these polynomials are holomorphic with respect to $w\in
D_1$. The analytic set defined by the polynomials $P_j$ contains $X$.
In particular, we have $P_j(x_j,\pi(x))=0$.  We consider the
case where $z\not=\pi(x)$, otherwise the lemma is clear.
We will show the existence of a
$z^{(i)}$ with good estimates on $z^{(i)}_j-x_j$.

Fix a constant $c>1$ large enough. There is an integer $2\leq l\leq
4m(k-p)+2$ such that $P_j(t,\pi(x))$ has no root $t$ with
$$(l-1)c\sqrt[m]{\|z-\pi(x)\|} < |t-x_j|\leq (l+1)c
\sqrt[m]{\|z-\pi(x)\|}$$
for every $p+1\leq j\leq k$. 
We call this {\it the security ring}.
For $\theta\in\R$ define
$$\xi_j:=lc\sqrt[m]{\|z-\pi(x)\|} e^{i\theta}+x_j$$
and
$$G_{j,c,\theta}(w):=c^{-m+1}\prod_{i=1}^m (\xi_j -
w_j^{(i)})=c^{-m+1}P_j(\xi_j,w).$$
Observe that $G_{j,c,\theta}(w)$ are Lipschitz with respect to $w$ in a
neighbourhood of $D$ uniformly with respect
to $(j,c,\theta)$. Using the choice of $l$, we have
$$|G_{j,c,\theta}(\pi(x))|=c^{-m+1}|P_j(\xi_j,\pi(x))|\geq c\|z-\pi(x)\|.$$
Hence, if $c$ is large enough, since the $G_{j,c,\theta}(w)$ are
uniformly Lipschitz, they do not vanish
on the ball $\widetilde D$ of center $\pi(x)$ and of radius
$2\|z-\pi(x)\|$.
Note that here we only need to consider the
  case where $z$ and $\pi(x)$ are close enough, and we have
  $\widetilde D\Subset
  D_1$.  
We denote by $\Sigma$ the boundary of the polydisc $H$ of
center $(x_{p+1},\ldots,x_k)\in\C^{k-p}$ and of radius
$lc\sqrt[m]{\|z-\pi(x)\|}$ : the $P_j(t,w)$ have no zero there when  $w\in
\widetilde D$. Then, $X$ does not intersect $\widetilde D\times
  \Sigma$. 
Since $z\in \widetilde D$ and $x\in X$, by continuity, there is a point $z^{(i)}$ satisfying
  $|z^{(i)}_j-x_j|\leq lc\sqrt[m]{\|z-\pi(x)\|}$. This gives the first
  assertion of the lemma. 

We now prove the second assertion. Fix a point $x$ in $\pi^{-1}(y)\cap
X$ and use the above construction. In the box $\widetilde D\times H$,
$X$ is a ramified covering over $\widetilde D$ of some degree $s\leq
m$. So we can write with an arbitrary order
$$\pi^{-1}(y)\cap X\cap \widetilde D\times H=\{y^{(1)},\ldots, y^{(s)}\}\quad \mbox{and} \quad
\pi^{-1}(z)\cap X\cap \widetilde D\times H =\{z^{(1)},\ldots,
z^{(s)}\}$$ 
with the desired estimates on $|y^{(i)}-z^{(i)}|$, since the diameter
of $\widetilde D\times H$ is controled by $\|y-z\|^{1/m}$. This gives a
partial correspondence between $\pi^{-1}(y)\cap X$
and $\pi^{-1}(z)\cap X$.

Choose another point $x'\in \pi^{-1}(y)\cap X$ outside $\widetilde
D\times H$ and repeat the construction in order to obtain a box $\widetilde D\times
H'$. We only replace the constant $c$ by $8[m(k-p)+1]c$. This
garantees that either $D\times H$ and  $D\times H'$ are disjoint or $\widetilde D\times H$ is
contained in $\widetilde D\times H'$ because of the security rings. In the last situation, we remove the 
box  $\widetilde D\times H$. Then, we repeat the construction for
points outside the boxes obtained so far. After less than $m$ steps, we
obtain a finite family of boxes which induces 
a complete correspondence between $\pi^{-1}(y)\cap X$
and $\pi^{-1}(z)\cap X$ satisfying the lemma.
\endproof

\begin{lemma} We have $\nu(u',x)< m^p\delta$ for every $u'\in\Gc'$ and $x\in D_1$.
\end{lemma}
\proof
Consider the functions $u'\in \Gc'$ and $u''$ as above, see (\ref{eq_psh_push}).
Let $y$ be a point in $\pi^{-1}(x)\cap X$ and $V$ a neighbourhood of
$y$ such that $\pi^{-1}(x)\cap X\cap V=\{y\}$. We can choose $V$ so
that $X\cap V$ is a ramified covering over $\pi(V)$. Let $l$ denote
the degree of this covering.
Consider the current $R:=\ddc(u[X])$ in $V$. In a
neighbourhood of $x$, $\ddc u'$ (which is equal to $\ddc u''$) is the
sum of the currents $\pi_*(R)$ for $y$ varying in $\pi^{-1}(x)\cap X$.
Since $\nu(R,y)<\delta$ and $l\leq m$, it is enough to prove that
$\nu(\pi_*(R),x)\leq l^{p-1}\nu(R,y)$. Assume that $y=0$ and $x=0$ in order to simplify
the notation. If $z=(z',z'')=(z_1,\ldots,z_p,z_{p+1},\ldots, z_k)$ 
denote the coordinates in $\C^k=\C^p\times\C^{k-p}$, then the mass of
$\pi_*(R)\wedge (\ddc\log\|z'\|)^{p-1}$ at $x=0$ is equal to
$\nu(\pi_*(R),0)$. It follows from 
the definition of $\pi_*$ that the mass of $R\wedge
(\ddc\log\|z'\|)^{p-1}$ at $y=0$ is also equal to $\nu(\pi_*(R),0)$.
Define $v:=\max(\log\|z'\|,l\log\|z''\|-M)$ with $M>0$ large
enough. Lemma \ref{lemma_lojasiewicz} applied to $X\cap V$ implies that $v=\log\|z'\|$
on $X\cap V$. Hence,
$R\wedge (\ddc\log\|z'\|)^{p-1}=R\wedge (\ddc v)^{p-1}$. Since $v\geq
l\log\|z\|-M'$, $M'>0$,  the comparison lemma in \cite{Demailly2} implies
that the mass of $R\wedge (\ddc v)^{p-1}$ at 0 is smaller than the
mass of $l^{p-1}R\wedge (\ddc \log\|z\|)^{p-1}$ at $0$ which is
equal to $l^{p-1}\nu(R,0)$. This completes the proof.
\endproof

\noindent
{\bf End of the proof of Proposition \ref{prop_psh_log}.} Now, we apply the case of smooth variety to $\Gc'$.  
If $0<\rho<1$ then $\sup_B u'\geq m^p\delta \log
\rho-\const$, where $B$ is the ball of center $\pi(a)$ and of radius $\rho$ in
$\C^p$. Let $B'$ be the connected component of $X\cap\pi^{-1}(B)$ which contains
$a$. This is a ramified covering over $B$. Since $u$ is negative, we
have $\sup_{B'} u\geq \sup_B u''\geq\sup_B u'$, see Proposition \ref{prop_psh_equi}(b). 
Lemma \ref{lemma_lojasiewicz} implies that $B'$ is contained in the
union of the balls of center in $\pi^{-1}(\pi(a))\cap X$ and of radius
$A\rho^{1/m}$, $A>0$. 
In this union, consider the connected component containing $a$. It has
diameter $\leq 2m A\rho^{1/m}$.
Hence, $B'$ is contained in the ball $B^X_a(r)$ of center $a$ and
of radius $r:=2mA\rho^{1/m}$ in $X$. We have
$$\sup_{B^X_a(r)} u\geq   m^p\delta \log \rho -\const \geq
m^{p+1}\delta\log r-\const$$
for $0<\rho<1$. This gives the estimate in the proposition with $c=m^{p+1}$.
\endproof

Consider the case where $X$ is an analytic subset of pure dimension
$p$ of $\P^k$. The
following proposition is a direct consequence of the last one.

\begin{proposition} \label{prop_T_psh_exp}
Let $\Gc\subset L^1(X)$ be a compact family of modulo $T$ wpsh functions
on $X$. Let $\delta> 0$ such that
$\nu_X(u,x)< \delta$ for $u\in \Gc$ and $x\in X$.
Then, there exist constants $c>0$ and $A>0$ such that
$$\sup_{B^X_a(r)}u\geq c\delta\log r -A \quad \mbox{ for } u\in\Gc,\  a\in
X \mbox{ and } 0<r<1.$$
Moreover, the constant $c$ is independent of $\Gc$ and of $\delta$.
\end{proposition}

The following result is a consequence of an inequality due to Demailly
and M{\'e}o \cite{Demailly2, Meo2}. It gives a bound for the volume of the
set where the Lelong numbers are large.

\begin{lemma} \label{lemma_demailly_meo}
Let $u$ be a modulo $T$ wpsh function on an analytic set  $X$  of pure
dimension $p$ in $\P^k$. Let $\beta\geq 0$ be a constant
and $q$ the dimension of  $\{\nu_X(u,x)>\beta\}$. Consider a finite
family of analytic sets $Z_r$, $1\leq r\leq s$, of pure dimension $q$ in $X$.
Assume that $\nu_X(u,x)\geq \nu_r$ for $x\in Z_r$ where
$(\nu_r)$ is a decreasing sequence such that $\nu_r\geq 2\beta$. 
Assume also that $\deg Z_r\geq d_r$ where the 
$d_r$'s are positive and satisfy $d_{r-1}\leq \frac{1}{2}d_r$. 
Then 
$$\sum_r d_r\nu_r^{p-q} \leq 2^{p-q+1}\deg(X)^{p-q}.$$
\end{lemma}
\proof
Define $R:=\ddc (u[X])+T\wedge [X]$, then $R$ is of bidimension $(p-1,p-1)$. Recall that
$\nu_X(u,x)=\nu(R,x)$. The mass of $R$ is equal to $\deg(X)$. 
Define $Z_1':=Z_1$ and for $r\geq 2$, $Z_r'$ the union of irreducible 
components of $Z_r$ which are not
components of $Z_1\cup\ldots\cup Z_{r-1}$.
So,  $Z_i'$ and $Z_r'$ have
no common component for $i\not=r$.
Let $d_r'$ denote the degree of $Z_r'$. We have 
$d_1'+\cdots+d_r'\geq d_r$ for $r\geq 1$. We also have $\nu(R,x)\geq \nu_r$
on $Z_r'$.
The inequality of Demailly-M{\'e}o
\cite{Demailly2, Meo2} implies that 
$$\sum_r (\deg Z_r') (\nu_r-\beta)^{p-q}\leq \|R\|^{p-q} = (\deg
X)^{p-q}.$$
Hence, since  $\beta\leq \nu_r/2$, 
$$\sum_r d_r' \nu_r^{p-q}\leq 2^{p-q} (\deg X)^{p-q}.$$
On the other hand, using the properties of $d_r$, $d_r'$, the fact that $(\nu_r)$ is
decreasing and the Abel's transform, we obtain
\begin{eqnarray*}
\sum_r d_r'\nu_r^{p-q} & = &
d_1'(\nu_1^{p-q}-\nu_2^{p-q})+(d_1'+d_2')(\nu_2^{p-q}-\nu_3^{p-q})
+\cdots+\\
& & + (d_1'+\cdots + d_{s-1}')(\nu_{s-1}^{p-q}-\nu_s^{p-q}) +
(d_1'+\cdots+d_s')\nu_s^{p-q} \\
&\geq & d_1(\nu_1^{p-q}-\nu_2^{p-q})+\cdots+d_{s-1}(\nu_{s-1}^{p-q}-\nu_s^{p-q}) +
d_s\nu_s^{p-q} \\
& \geq & \frac{1}{2} d_1\nu_1^{p-q}+\cdots+\frac{1}{2} d_s\nu_s^{p-q}.
\end{eqnarray*}
This proves the lemma.
\endproof


\section{Asymptotic contraction}

In this section, we study the speed of contraction of $f^n$.
More precisely, we want to estimate the size of the largest ball
contained in the image of a fixed ball by $f^n$. Our main result
is the following theorem where the balls in $X$ are defined in Section
\ref{section_lelong}.

\begin{theorem} \label{th_contraction}
Let $f$ be a holomorphic endomorphism of algebraic degree $d\geq
2$ of $\P^k$ and $X$ an analytic subset of pure dimension $p$, $1\leq
p\leq k$, invariant
by $f$, i.e. $f(X)=X$. There exists a  constant $c>0$  such that if
$B$ is a ball of radius $r$ in $X$ with $0<r<1$, then for every $n\geq 0$, $f^n(B)$ contains a
ball in $X$ of radius $\exp(- cr^{-2p}d^n)$.
\end{theorem}

\begin{corollary} \label{cor_contraction}
Let $f$ be a holomorphic endomorphism of algebraic degree $d\geq
2$ of $\P^k$. There exists a  constant $c>0$  such that if
$B$ is a ball of radius $r$ in $\P^k$ with $0<r<1$, then $f^n(B)$ contains a
ball of radius $\exp(- cr^{-2k}d^n)$ for every $n\geq 0$.
\end{corollary}

Let $H$ be a hypersuface in $\P^k$ which does not contain any
irreducible component of $X$ such that the restriction of $f$ to $X\setminus H$ is
of maximal rank at every point. 
We choose $H$ containing $\sing(X)\cup f^{-1}(\sing(X))$.
If $\delta$ is the degree of $H$, there is a
negative function $u$ on $\P^k$ psh modulo $T$ such that $\ddc u= \delta^{-1}[H]- T$.

\begin{lemma} \label{lemma_image_ball}
There are positive constants  $c_1$ and $c_2$ such that if $B$ is a ball of center $a$ and
of radius $0<r<1$ in $X$ then $f(B)$ contains the ball 
of center $f(a)$ and of radius $c_1r\exp(c_2u(a))$ in $X$.
Moreover, if $u(a)\not=-\infty$ then the differential at $f(a)$ of $f^{-1}$
restricted to $X$ satisfies $\|Df_{|X}^{-1}(f(a))\|\leq c_1^{-1}\exp(-c_2u(a))$.
\end{lemma}
\proof
The constants $c_i$ that we use here are independent of $a$ and $r$.
We only have to consider the case where  $u(a)\not=-\infty$.
Observe that when $c_1$ is small and $c_2$ is large enough, the ball
of center $f(a)$ and of radius $c_1r\exp(c_2u(a))$ does not intersect $\sing(X)$.
Let $\pi:\widetilde X\rightarrow X\subset\P^k$ be a desingularization of $X$
and $A:=\|\pi\|_{\Cc^1}$. If $\pi(\widetilde a)=a$, and if
$\widetilde B$ is the ball of center $\widetilde a$ and of radius
$\widetilde r:=A^{-1}r$ then $\pi(\widetilde B)$ is contained in the
ball $B$. 
Define  $h:=f\circ\pi$, $\widetilde u:=u\circ\pi$
and $\widetilde T:=\pi^*(T)$. Since $T$ has continuous local
potentials, so does $\pi^*(T)$.

The current $\pi^*[ H]$ is supported in $\pi^{-1}(H)$ and satisfies
$\ddc \widetilde u = \delta^{-1}\pi^*[H]-\pi^*(T)$. Since $\widetilde u=-\infty$
exactly on $\pi^{-1}(H)$ and since $\pi^*(T)$ has continuous local
potentials, the support of $\pi^*[H]$ is exactly $\pi^{-1}(H)$. So,
$\pi^*[H]$ is a combination with strictly positive coefficients of the
currents of
integration on irreducible components of $\pi^{-1}(H)$. Observe that $h$ is of
maximal rank outside $\pi^{-1}(H)$.
It is enough to prove that $h(\widetilde B)$ contains the ball of
center $h(\widetilde a)$
and of radius $c_1r\exp(c_2\widetilde u (\widetilde a))$ in $X$.

We can assume that $r$ is small and work in the local setting. We
use holomorphic coordinates $x=(x_1,\ldots,x_p)$ of $\widetilde X$ and
$y=(y_1,\ldots,y_k)$ of $\P^k$ in small
neighbourhoods $W$ and $U$ of $\widetilde a$ and $a$ respectively. 
Write $h=(h_1,\ldots,h_k)$ and consider 
a holomorphic function $\varphi$ on $W$ such that
$\varphi^{-1}(0)=\pi^{-1}(H)\cap W$. Then, $\delta^{-1}\pi^*[H]\geq
\epsilon\ddc\log|\varphi|$ with $\epsilon>0$ small enough.
We have $\ddc(\widetilde u - \epsilon\log|\varphi|)\geq -\widetilde T$.
 It follows that  
$\widetilde u-\epsilon \log|\varphi|$ is a difference of a psh function and a
potential of $\widetilde T$. Since
$\widetilde T$ has local continuous potentials, $\widetilde
u-\epsilon\log|\varphi|$ is bounded from above.
Up to multiplying $\varphi$ by a constant, we can assume
that $\epsilon \log|\varphi|\geq \widetilde u$.

If $J\subset \{1,\ldots,k\}$ is a multi-index of length $p$, denote by
$M_J$ the matrix $(\partial h_j/\partial x_i)$ with $1\leq i\leq p$
and $j\in J$. 
Since $h$ is of maximal rank outside $\pi^{-1}(H)$, 
the zero set of $\sum_J|\det M_J|^2$ is contained in $\{\varphi=0\}$.
The Lojasiewicz's inequality \cite{Tougeron} implies that $\sum_J
|\det M_J|^2 \geq c_3|\varphi|^{c_4}$ for some constants $c_3>0$ and $c_4>0$.
Up to a
permutation of the coordinates $y$, we can assume that $|\det
M(\widetilde a)|\geq c_5|\varphi(\widetilde a)|^{c_2\epsilon/2}\geq
c_5\exp(c_2\widetilde u(\widetilde a)/2)$ where $c_2$, $c_5$ are
positive constants and $M$ is
the matrix $(\partial h_j/\partial x_i)_{1\leq i,j\leq p}$.
Define $h':=(h_1,\ldots,h_p)$. The precise version of the implicit
function theorem  \cite[p.106]{Tougeron} implies that $h'$ defines a
bijection from an open subset of $\widetilde B$ to a ball of center
$h'(\widetilde a)$ and of radius $c_6\widetilde r|\det M(\widetilde
a)|^2$, $c_6>0$. This proves the first assertion in the lemma. For the
second one, we have $\|Dh'^{-1}\|\lesssim |\det M(\widetilde
a)|^{-1}$ at $h'(\widetilde a)$ which gives the result.
\endproof

\noindent {\bf Proof of Theorem \ref{th_contraction}.} 
By Corollary \ref{cor_psh_convergence_T}, 
the sequence of functions 
$(d^{-n}u\circ f^n)$ is bounded
in $L^1(X)$. Since 
$$d^{-n}(u+u\circ f+\cdots+u\circ
f^{n-1})=\sum_{i=0}^{n-1}d^{-(n-i)}(d^{-i}u\circ f^i),$$ 
the $L^1(X)$-norm of $d^{-n}(u+u\circ f+\cdots+u\circ
f^{n-1})$ is bounded by a constant $c'>0$ independent of $n$.
If $A>0$ is
a constant large enough, the set of points $x\in X$ satisfying
$$u(x)+u\circ f(x)+\cdots+u\circ
f^{n-1}(x)\leq -Ar^{-2p}d^n$$
has Lebesgue measure $\leq c'A^{-1}r^{2p}$. 
By a theorem of Lelong \cite{Lelong, Demailly2}, the volume of a ball
of radius $r/2$ in $X$ is $\geq c''r^{2p}$, $c''>0$.
Therefore, since $A$ is
large, there  is a point $b\in X$,
depending on $n$,  such that $|b-a|\leq r/2$ and 
\begin{equation}\label{eq_radius}
u(b)+u\circ f(b)+\cdots+u\circ f^{n-1}(b)\geq -Ar^{-2p}d^n.
\end{equation}
Lemma \ref{lemma_image_ball} applied inductively to
balls centered at $f^i(b)$  implies that $f^n(B)$ contains the
ball of center $f^n(b)$ of radius
$$\frac{1}{2}c_1^nr \exp\big(c_2u(b)+\cdots+c_2u(f^{n-1}(b))\big).$$
We obtain the result using (\ref{eq_radius}) and the estimate $\frac{1}{2}c_1^nr\geq
\exp(-c_3r^{-2p}d^n)$ for $0<r<1$, where $c_3>0$ is a constant.
\hfill $\square$

\begin{remark} \rm
With the same argument we also get the following. Let $B_x$ denote the
ball of center $x$ and of radius $0<r<1$ in $X$. Let $r_n(x)$ be the
maximal radius of the ball centered at $f^n(x)$ and contained in
$f^n(B_x)$. Then, there is a constant $A>0$ such that 
$$\frac{\log r_n(x)}{d^n}\geq -\frac{A(n+1)-\log r}{d^n}+\frac{c_2}{d^n}\sum_{i=0}^{n-1}
u\circ f^i(x).$$
Consequently, there is a constant $c>0$ such that
$$\int_X \frac{\log r_n(x)}{d^n} \omega^p\geq -c+\frac{\log
  r}{d^n}\deg X.$$
We can also replace $\omega^p$ by any PB measure on $X$, i.e. a
measure such that modulo $T$ wpsh functions are integrable, see
\cite{DinhSibony4}. 
\end{remark}

In the following result, we use the Lebesgue measure $\volume_X$ on $X$ induced by the
Fubini-Study form restricted to $X$.

\begin{theorem} \label{th_volume}
Let $f$ and $X$ be as in Theorem \ref{th_contraction}.  Let $Z$ be a
  Borel set in $X$ and $n\geq 0$. Then there is a Borel set $Z_n\subset Z$ with
  $\volume_X(Z_n)\geq {1\over 2}\volume_X(Z)$ such that the
  restriction $f_{|X}^n$ of $f^n$ to $X$  defines a
  locally bi-Lipschitz map from $Z_n$ to $f^n(Z_n)$. Moreover, the
  differential of the inverse map $f_{|X}^{-n}$ satisfies
  $\|Df_{|X}^{-n}\|\leq \exp(c\volume(Z)^{-1}d^n)$ on $f^n(Z_n)$ with
  a constant $c>0$ independent of $n$ and  $Z$. In particular, we have
    $\volume_X(f^n(Z))\geq \exp(-c'\volume_X(Z)^{-1}d^{n})$ for some constant
      $c'>0$ independent of $n$ and  $Z$.
\end{theorem}
\proof
As in (\ref{eq_radius}), there is a subset $Z_n$ of $Z$ with
$\volume_X(Z_n)\geq {1\over 2}\volume_X(Z)$ such that 
$$u(b)+u\circ f(b)+\cdots+u\circ f^{n-1}(b)\geq -A\volume_X(Z)^{-1}d^n$$
for $b\in Z_n$, where $A>0$ is a fixed constant large enough. 
In particular, we have $u\circ f^i(b)\not=-\infty$ for $i\leq n-1$. It
follows from the definition of $u$ that $f_{|X}^n$ defines a bijection
between a neighbourhood of $b$ and a neighbourhood of $f^n(b)$ in $X$. Hence,
$f_{|X}^n:Z_n\rightarrow f^n(Z_n)$
is locally bi-Lipschitz. Applying Lemma \ref{lemma_image_ball} inductively gives the
estimate on  $\|Df_{|X}^{-n}\|$ at $f^n(b)$. 

Since the fibers of $f^n$ contains at most $d^{kn}$ points, the
estimate on  $\|Df_{|X}^{-n}\|$ implies
$$\volume_X(f^n(Z))\geq \volume_X(f^n(Z_n))\gtrsim d^{-kn}\volume_X(Z)\exp(-c\volume_X(Z)^{-1}d^{n})^{2p}.$$
The last assertion in the theorem follows.
\endproof

\begin{remark}\rm
It is not difficult to extend Theorems \ref{th_contraction} and
\ref{th_volume} to the case
of meromorphic maps or correspondences on compact K{\"a}hler manifolds. We can use the continuity
of $f^*$ on the space $\DSH$ in order to estimate the $L^1$-norm
of $u\circ f^n$ for $u\in\DSH$, see \cite{DinhSibony4}. The volume
estimate in Theorem \ref{th_volume} for meromorphic maps on
smooth manifolds was obtained in \cite{Guedj2}, see also 
\cite{FornaessSibony1, FavreJonsson, Guedj} for earlier versions. 
\end{remark}

Let $\Gc$ be a compact family of modulo $T$ wpsh functions on $X$. 
Let $\Hc_n$ denote the family of $T$ wpsh functions which
are equal almost everywhere to $d^{-n} u\circ f^n$, $u\in\Gc$. 
Define 
$$\nu_{n}:=\sup\{\nu_X(u,a),\ u\in\Hc_{n},\ a\in X\}.$$
We have the following result.

\begin{proposition} \label{prop_limit_lelong}
Assume that $\inf \nu_{n}=0$. Then,  $d^{-n}
u_{n}\circ f^{n}\rightarrow 0$ in $L^1(X)$ for
all $u_{n}\in\Gc$. In particular, the hypothesis is satisfied when there is an
increasing sequence $(n_i)$
such that $d^{-n_i} u_{n_i}\circ f^{n_i}$ converge to $0$ in $L^1(X)$ for
all $u_{n_i}\in\Gc$.
\end{proposition}
\proof
Consider
a sequence $(d^{-n_i} u_{n_i}\circ f^{n_i})$ converging in $L^1(X)$ to a modulo
$T$ wpsh function $u$. Corollary \ref{cor_psh_convergence_T} implies that $u\leq 0$. 
We want to prove that $u=0$. If not, since $u$ is upper
semi-continuous, there is a constant $\alpha>0$ such
that $u\leq -2\alpha$ on some ball $B$ of radius $0<r<1$ in $X$. 
By Lemmas \ref{lemma_hartogs_T}
and \ref{lemma_t_psh_interate},
for $i$ large enough, we have $d^{-n_i}u_{n_i}\circ f^{n_i}\leq
-\alpha$ almost everywhere on $B$. 

Fix $\delta>0$ small enough 
and $m$ such that $\nu_m<\delta$. Consider only the $n_i$ larger than $m$.
Then,
$d^{-m}u_{n_i}\circ f^m\leq
-d^{n_i-m}\alpha$ almost everywhere on $f^{n_i-m}(B)$. By
Theorem \ref{th_contraction},  $f^{n_i-m}(B)$ contains a ball $B_i$ of radius
$\exp(-cr^{-2p} d^{n_i-m})$ in $X$ with $c>1$. If $v_i\in\Hc_m$ is equal
almost everywhere to $d^{-m} u_{n_i}\circ f^m$, then $v_i\leq
-d^{n_i-m}\alpha$ almost everywhere on $B_i$. It follows from
Proposition \ref{prop_equi_T_psh}(b) that this inequality holds 
everywhere on $B_i$.
By Proposition \ref{prop_T_psh_exp}, there is a constant $c'>0$
independent of $\Gc$, $r$, $\delta$, $m$, and a constant $A>0$ such that
$$-c'\delta r^{-2p}d^{n_i-m}-A\leq \sup_{B_i} v_i \leq -d^{n_i-m}\alpha.$$
This is a contradiction if $\delta$ is chosen small enough and if $n_i$
is large enough.

Assume now that $d^{-n_i} u_{n_i}\circ f^{n_i}$ converge to $0$ in $L^1(X)$ for
all $u_{n_i}\in\Gc$.  Then, for every $\epsilon>0$, we
have $\nu(u,a)<\epsilon$ for $u\in  \Hc_{n_i}$, $a\in X$ and
for $i$ large enough. Therefore, $\inf\nu_n=0$. Here, we use that if
positive closed currents $R_n$ converge to $R$ and $a_n\rightarrow a$
then $\limsup\nu(R_n,a_n)\leq\nu(R,a)$. 
\endproof

\begin{corollary} 
Let $\Fc$ be a family of positive closed $(1,1)$-currents of mass
  $1$ on $\P^k$. Assume that there is an increasing sequence of integers
  $(n_i)$ such that $d^{-n_i} (f^{n_i})^*(S_{n_i})\rightarrow T$ for all
  $S_{n_i}\in\Fc$. Then, $d^{-n} (f^n)^*(S_n)\rightarrow T$ for all $S_n\in\Fc$.
\end{corollary}
\proof
Observe that the hypothesis implies that $d^{-n_i} (f^{n_i})^*(S_{n_i})\rightarrow T$ for all
  $S_{n_i}\in\overline\Fc$. So, we can replace $\Fc$ by $\overline\Fc$ and
  assume that $\Fc$ is compact.
To each current $S\in\Fc$ we associate a modulo $T$ psh function $u$
on $\P^k$ such
that $\ddc u=S-T$. Subtracting from $u$ some constant allows us to
have $\max_{\P^k} u=0$. 
Proposition
\ref{prop_T_psh_compactness} and Lemma \ref{lemma_hartogs_T}
imply that the family $\Gc$ of these functions $u$ is compact.
The hypothesis and
Corollary \ref{cor_psh_convergence_T} imply that $d^{-n_i} u_{n_i}\circ f^{n_i}\rightarrow 0$
for $u_{n_i}\in\Gc$. Proposition \ref{prop_limit_lelong}
gives the result.
\endproof

\begin{corollary} \label{cor_lelong_zero}
Let $\Fc$ be a compact family of positive closed $(1,1)$-currents of
mass $1$ on $\P^k$. Assume that for any $S\in\Fc$, the Lelong number of $S$ vanishes at
every point out of $\supp(\mu)$. Then, $d^{-n}(f^n)^*(S_n)\rightarrow T$
for any sequence $(S_n)\subset \Fc$. 
\end{corollary}
\proof
Let $\Gc$ and $\Hc_n$ be defined as above.
Define also
$$\nu_n':=\sup\{\nu_X(u,a),\ u\in\Hc_n,\ a\in \supp(\mu)\}$$
and
$$\nu_n'':=\sup\{\nu_X(u,a),\ u\in\Hc_n,\ a\not\in \supp(\mu)\}.$$
Corollary \ref{cor_psh_convergence_T} implies that $\lim\nu_n'=0$.
On the other hand, by hypothesis, $\nu_0''=0$. Since $\P^k\setminus\supp(\mu)$ is totally
invariant, Proposition \ref{prop_lelong_compare}, applied to $X=\P^k$, implies that $\nu_n''=0$ for every $n$.
Hence, $\nu_n=\nu_n'$ and $\nu_n\rightarrow 0$.
We apply Proposition \ref{prop_limit_lelong} in order to conclude.
Note that the corollary still holds if we only assume
that $\inf \nu_n''=0$.
\endproof

We prove as in Proposition \ref{prop_limit_lelong} the following lemma.

\begin{lemma} \label{lemma_T_psh_lelong_zero}
Let $(u_{n_i})$ be a sequence of modulo $T$ wpsh functions on
$X$, bounded in $L^1(X)$. Assume that $d^{-n_i}u_{n_i}\circ f^{n_i}$ converge to a modulo
$T$ wpsh function $v$. Assume also that for every $\delta>0$, there is
a subsequence $(u_{m_i})\subset (u_{n_i})$ converging to a modulo $T$
wpsh function $w$ with $\nu_X(w,a)<\delta$ at every point $a\in X$. Then, $v=0$.
\end{lemma}
\proof
Corollary \ref{cor_psh_convergence_T} implies that $v\leq 0$. Assume that
$v\not=0$. Then, since $v$ is upper semi-continuous, there is a
constant $\alpha>0$ such that
$v<-2\alpha$ on a ball of radius $0<r<1$ on $X$.
As in Proposition  \ref{prop_limit_lelong}, for $i$ large enough we have
$u_{n_i}<-d^{n_i}\alpha$ on a ball $B_{n_i}$ of radius
$\exp(-cr^{-2p}d^{n_i})$ in $X$ with $c>1$.

Fix $\delta>0$ small enough, and $(u_{m_i})$ and $w$ as above. 
The property of $w$ implies that if $s$ is an integer
large enough, we have $\nu_X(u_{m_i},a)<\delta$ for every $a\in
X$ and for $i\geq s$. By Proposition \ref{prop_T_psh_exp} applied to the
compact family $\{u_{m_i},\ i\geq s\}\cup\{w\}$, there is a constant $c'>0$
independent of $\delta$, $r$ and a constant $A>0$ such that
$$-c'\delta r^{-2p}d^{m_i}-A\leq \sup_{B_{m_i}} u_{m_i} \leq
-d^{m_i}\alpha\quad \mbox{for } i\geq s.$$
This is a contradiction 
for $m_i$ large enough, since $\delta$ is chosen small.
\endproof


\section{Exceptional sets} \label{section_exceptional}

Let $X$ be an analytic subset of pure dimension $p$ in $\P^k$
invariant by $f$, i.e. $f(X)=X$. Let
$g:X\rightarrow X$ denote the restriction of $f$ to $X$. 
We will follow the idea of \cite{DinhSibony1} in order to define and study
the exceptional analytic subset $\Ec_X$ of $X$ which is
totally invariant by $g$, see also \cite{Dinh1, Dinh2}.
The following result can be deduced from Section 3.4
in \cite{DinhSibony1}.

\begin{theorem} \label{th_exceptional}
There is a (possibly empty) proper analytic subset $\Ec_X$ of $X$
which is totally invariant by $g$ and is maximal in the following
sense. If $E$ is an analytic subset of dimension $<p$ of $X$ such
that $g^{-s}(E)\subset E$ for some $s\geq 1$, then $E\subset \Ec_X$.
In particular, there is a maximal proper analytic subset $\Ec_{\P^k}$ of $\P^k$
which is totally invariant by $f$.
\end{theorem}

We will need some precise properties of $\Ec_X$. So, for the reader's
convenience, we recall here the
construction of $\Ec_X$ and the proof of the previous theorem since the emphasis
in \cite{DinhSibony1} is on polynomial-like maps.
Observe that $g$ permutes the irreducible components of $X$. Let $m\geq 1$ be an
integer such that $g^m$ fixes the components of $X$. 

\begin{lemma} \label{lemma_topological}
The topological degree of  $g^m$ is equal to 
$d^{mp}$, that is, $g^m:X\rightarrow X$ defines a ramified covering of
degree $d^{mp}$. In particular, for every $x\in X$, $g^{-m}(x)$ contains at
most $d^{mp}$ points and there is a hypersurface $Y$ of $X$ containing
$\sing(X)\cup g^m(\sing(X))$ such that for $x\in X\setminus Y$, $g^{-m}(x)$
contains exactly $d^{mp}$ points.
\end{lemma}
\proof
We can work
with each component. So, we can assume that $X$ is irreducible. It follows
that $g^m$ defines a ramified covering. We want to
prove that the degree $\delta$ of this covering is equal to $d^{mp}$.
Consider the positive measure
$(f^m)^*(\omega^p)\wedge[X]$. Its mass is equal to 
$d^{mp}\deg(X)$ since $(f^m)^*(\omega^p)$ is cohomologous to
$d^{mp}\omega^p$. The operator $(f^m)_*$
preserves the mass of positive measures. We also have
$(f^m)_*[X]=\delta[X]$. Hence,
\begin{eqnarray*}
d^{mp}\deg(X) & = & \|(f^m)^*(\omega^p)\wedge [X]\| =
\|(f^m)_*((f^m)^*(\omega^p)\wedge [X])\| \\
& = &
\|\omega^p\wedge (f^m)_*[X]\|  
=  \delta\|\omega^p\wedge [X]\|=\delta\deg(X).
\end{eqnarray*}
Therefore, $\delta=d^{mp}$.  So, we can take for $Y$, a hypersurface
containing the ramification values of $f^m$ and $\sing(X)\cup g^m(\sing(X))$.
\endproof

Let $Y$ be as above. Observe that if $g^m(x)\not\in Y$ then $g^m$ defines a biholomorphic
map between a neighbourhood of $x$ and a neighbourhood of $g^m(x)$ in
$X$. 
Let $[Y]$ denote the $(k-p+1,k-p+1)$-current of integration
on $Y$ in $\P^k$. Since $(f^{mn})_*[Y]$ is a positive closed
$(k-p+1,k-p+1)$-current of mass $d^{mn(p-1)}\deg(Y)$, we can define
the following ramification current
$$R=\sum_{n\geq 0} R_n:=\sum_{n\geq 0} d^{-mnp} (f^{mn})_*[Y].$$  
By a theorem of Siu \cite{Siu,Demailly2}, for $c>0$, the level set
$E_c:=\{\nu(R,x)\geq c\}$ of the Lelong number
is an analytic set of dimension $\leq p-1$ contained in $X$. Observe
that $E_1$  contains $Y$. We will see that $R$ is the obstruction for
the construction of ``regular'' orbits.

For any point $x\in X$ let $\lambda'_n(x)$ denote the number of  distinct orbits
$$x_{-n},x_{-n+1},\ldots,x_{-1}, x_0$$
such that $g^m(x_{-i-1})=x_{-i}$, $x_0=x$ and $x_{-i}\in X\setminus Y$ for
$0\leq i\leq n-1$. These are the ``good" orbits. Define $\lambda_n:=d^{-mpn}\lambda'_n$.
The function $\lambda_n$ is lower semi-continuous with respect to
the Zariski topology on $X$. Moreover, by Lemma
\ref{lemma_topological}, 
we have $0\leq \lambda_n\leq 1$ and $\lambda_n=1$ 
out of the analytic set $\cup_{i=0}^{n-1}g^{mi}(Y)$. The sequence 
$(\lambda_n)$ decreases to a function $\lambda$, which represents the
asymptotic proportion of orbits in $X\setminus Y$.

\begin{lemma} \label{lemma_exceptional}
There is a constant $\gamma>0$ such that
  $\lambda\geq\gamma$ on   $X\setminus E_1$.
\end{lemma}
\proof
We deduce from the Siu's theorem,  the existence of a constant
$0<\gamma<1$ satisfying $\{\nu(R,x)> 1-\gamma\}=E_1$. Consider a
point $x\in X\setminus E_1$. We have $x\not\in Y$. Define $\nu_n:=\nu(R_n,x)$. We have $\sum
\nu_n\leq 1-\gamma$. 
Since $E_1$ contains $Y$, $\nu_0=0$ and $F_1:=g^{-m}(x)$
contains exactly $d^{mp}$ points. The definition of $\nu_1$ implies
that 
$g^{-m}(x)$ contains at most $\nu_1d^{mp}$ points in $Y$. 
Then 
$$\#g^{-m}(F_1\setminus Y)=d^{mp}\#(F_1\setminus Y) \geq
(1-\nu_1)d^{2mp}.$$
Define $F_2:=g^{-m}(F_1\setminus Y)$. The definition of $\nu_2$
implies that $F_2$ contains at most 
$\nu_2 d^{2mp}$ points in $Y$. Hence, $F_3:=g^{-m}(F_2\setminus Y)$
contains at least $(1-\nu_1-\nu_2)d^{3mp}$ points. In the same way, we
define $F_4$, $\ldots$, $F_n$ with $\#F_n\geq (1-\sum
\nu_i)d^{mpn}$. Hence, for every $n$ we get the following estimate: 
$$\lambda_n(x)\geq
d^{-mpn}\#F_n\geq 1-\sum\nu_i\geq \gamma.$$ 
This proves the lemma.
\endproof

\medskip

\noindent
{\bf End of the proof of Theorem \ref{th_exceptional}.}
Let $\Ec_X^n$ denote the set of  $x\in X$ such that
$g^{-ml}(x)\subset E_1$ for $0\leq l\leq n$ and define
$\Ec_X:=\cap_{n\geq 0 }\Ec_X^n$.
Then, $(\Ec_X^n)$ is a decreasing sequence of analytic subsets of $E_1$.
It should be stationary. So, there is $n_0\geq 0$ such that
$\Ec^n_X=\Ec_X$ for $n\geq n_0$.

By definition, $\Ec_X$ is the set of $x\in X$ such that
$g^{-mn}(x)\subset E_1$ for every $n\geq 0$. Hence,
$g^{-m}(\Ec_X)\subset \Ec_X$. It follows that the sequence of analytic
sets $g^{-mn}(\Ec_X)$ is decreasing and there is $n\geq 0$ such that
$g^{-m(n+1)}(\Ec_X)=g^{-mn} (\Ec_X)$. Since $g^{mn}$  is surjective,
we deduce that $g^{-m}(\Ec_X)=\Ec_X$ and hence $\Ec_X=g^m(\Ec_X)$.

Assume as in the theorem that $E$ is analytic with $g^{-s}(E)\subset
E$. Define $E':=g^{-s+1}(E)\cup\ldots\cup E$. We have $g^{-1}(E')\subset
E'$ which implies $g^{-n-1}(E') \subset g^{-n}(E')$ for every $n\geq
0$. Hence, $g^{-n-1}(E')=g^{-n}(E')$ for $n$ large enough. This and
the surjectivity of $g$ imply that
$g^{-1}(E')=g(E')=E'$. 
By Lemma \ref{lemma_topological}, the topological degree of
$(g^{m'})_{|E'}$ is at most $d^{m'(p-1)}$ for some integer $m'\geq 1$.
This, the identity $g^{-1}(E')=g(E')=E'$ together with Lemma \ref{lemma_exceptional}
imply that $E'\subset E_1$. Hence,
$E'\subset \Ec_X$ and $E\subset \Ec_X$. 

Define $\Ec_X':=g^{-m+1}(\Ec_X)\cup\ldots\cup \Ec_X$. We have
$g^{-1}(\Ec_X')=g(\Ec_X')=\Ec_X'$. Applying the previous assertion to
$E:=\Ec'_X$ yields $\Ec'_X\subset \Ec_X$. Therefore,
$\Ec_X'=\Ec_X$ and $g^{-1}(\Ec_X)=g(\Ec_X)=\Ec_X$. 
\hfill $\square$

\begin{remark} \label{rk_exception}
\rm
The maximality of $\Ec_X$ in Theorem \ref{th_exceptional} implies that it does not
depend on the choice of $m$ and of the analytic set $Y$ satisfying Lemma
\ref{lemma_topological}. Moreover, $\Ec_X$ is also the exceptional set
associated to $g^n$ for every $n\geq 1$. An analytic set, totally
invariant by $g^n$, is not necessarily totally invariant by $g$, but
it is a union of components of such sets. We deduce from our
construction that $\Ec_{\P^k}$ depends algebraically on $f$.
\end{remark}

\begin{corollary} \label{cor_excep_finite}
There are only finitely many analytic subsets of $X$ which are totally
invariant by $g$. In particular, there is only a finite number of analytic
subsets of $\P^k$ which are totally invariant by $f$.
\end{corollary}
\proof
We only have to consider totally analytic sets $E$ of pure dimension
$q$. The proof is by induction on the dimension $p$ of $X$. Assume
that the corollary is true for $X$ of dimension $\leq p-1$ and
consider the case of dimension $p$. If $q=p$ then $E$ is a union of
components of $X$. There is only a finite number of such analytic
sets. If $q< p$, by Theorem \ref{th_exceptional}, $E$ is contained in
$\Ec_X$. Applying the hypothesis of induction to the restriction of
$f$ to $\Ec_X$ gives the result.
\endproof

We now give another characterization of $\E_X$.  Recall
that $\mu_X:=T^p\wedge [X]$. This is a positive measure of mass
$s:=\deg X$. The invariance of $T$ implies that $\mu$
is totally invariant by $g^m$, that is,
$(g^m)^*(\mu)=d^{pm}\mu$. Since $g^m$ fixes the components of $X$, we
can apply the to each component following result
where the
second assertion was proved by the authors in \cite{DinhSibony1}.

\begin{theorem} \label{th_conjecture_k}
Assume that $X$ is irreducible.
Let $\delta_a$ denote the Dirac mass at a point $a\in X$. Then
$d^{-pmn}(g^{mn})^*(\delta_a)$ converge to $s^{-1} \mu_X$ if and only if
$a$ is out of $\Ec_X$. In particular, if $a$ is a point in $\P^k$ then 
$d^{-kn}(f^n)^*(\delta_a)$ converge to $\mu$ if and only if
$a$ is out of $\Ec_{\P^k}$.
\end{theorem}

Since $T$ has continuous local potentials, $\mu_X$ has no
mass on proper analytic subsets of $X$. It follows that if
$a\in\Ec_X$, any limit value of $d^{-pmn}(g^{mn})^*(\delta_a)$ has support
in $\Ec_X$ and is singular with respect to $\mu_X$. Consider a point
$a$ in $X\setminus \Ec_X$. We only have to check the convergence to $s^{-1}\mu_X$.
Forn\ae ss and the second author proved this convergence for
$X=\P^k$ and for $a$ outside a pluripolar set 
\cite{FornaessSibony1}. Briend and Duval extended
this result to $a$ outside the orbit of the critical set of $f$
\cite{BriendDuval}. They also proposed a geometrical approach in order
to prove this
property for $a$ outside an analytic set but there is a problem with
the counting of multiplicity in their lemma in \cite[p.149]{BriendDuval}.

Briend-Duval result can be extended to our situation: for $a$ outside
the orbit of $Y$ we have $d^{-pmn} (g^{mn})^*(\delta_a)\rightarrow s^{-1} \mu_X$.
We recall the following proposition, see \cite{BriendDuval} and also
\cite{DinhSibony1, Dinh1, Dinh2} for more general cases, in
particular, for non-projective manifolds. 

\begin{proposition} \label{prop_Briend_Duval}
For any $\epsilon>0$, there is an integer
  $n_\epsilon\geq 0$ such that if
$a$ is out of the analytic set $Y_\epsilon:=Y\cup g^m(Y)\cup\ldots\cup
g^{mn_\epsilon}(Y)$, then
any limit value  $\nu$ of $d^{-pmn}(g^{mn})^*(\delta_a)$  satisfies 
$\|\nu-s^{-1}\mu_X\|\leq \epsilon$, where $\|\cdot\|$ denotes the mass
of measure.
\end{proposition}

Observe that if $n\geq r\geq 0$ then
$$d^{-pmn} (g^{mn})^*(\delta_a) = d^{-pmr} \sum_{b\in g^{-mr}(a)}
d^{-pm(n-r)} (g^{m(n-r)})^*(\delta_b),$$
where the points in $g^{-mr}(a)$ are counted with multiplicities.
So, if a point $a$ does not satisfies the conclusion of Proposition
\ref{prop_Briend_Duval} 
then it admits many preimages in $Y_\epsilon$. We quantify now this property.

Let $N_n(a)$ denote the number of orbits of $g^m$
$$\Oc=\{a_{-n},\ldots,a_{-1},a_0\}$$
with $g^m(a_{-i-1})=a_{-i}$ and $a_0=a$ such that $a_{-i}\in Y_\epsilon$
for every $i$. Here, the orbits are counted with multiplicities. So,
$N_n(a)$ is the number of
negative orbits of order $n$ of $a$ which stay in $Y_\epsilon$. Observe
that the sequence of functions $\tau_n:=d^{-pmn} N_n$ decreases to some
function $\tau$. Since $\tau_n$ are upper semi-continuous with respect
to the Zariski topology and $0\leq \tau_n\leq 1$, the function $\tau$
satisfies the same properties. 
Observe that $\tau(a)$ is the probability that an infinite negative
orbit of $a$ stays in $Y_\epsilon$.
The following proposition gives also a
characterization of $\Ec_X$.

\begin{proposition} The function $\tau$ is the characteristic function
  of $\Ec_X$, that is, $\tau=1$ on $\Ec_X$ and $\tau=0$ on $X\setminus
  \Ec_X$.
\end{proposition}
\proof
Since $\Ec_X\subset Y_\epsilon$ and $\Ec_X$ is totally invariant by
$g$, we have $\Ec_X\subset \{\tau=1\}$. Let $\theta\geq 0$ denote the
maximal value of $\tau$ on $X\setminus \Ec_X$. This value exists since
$\tau$ is upper semi-continuous with respect to the Zariski topology
(indeed, it is enough to consider the algebraic subset
$\{\tau\geq\theta_0\}$ of $X$ which decreases when $\theta_0$ increases).
We have to check that $\theta=0$. Assume in order to obtain a
contradiction that $\theta>0$. Since $\tau\leq 1$, we always have
$\theta\leq 1$. 
Consider the non-empty analytic set
$E:=\tau^{-1}(\theta)\setminus\Ec_X$ in $Y_\epsilon$. Let $a'$ be a point in $E$. Since
$\Ec_X$ is totally invariant, we have
$g^{-m}(a')\cap\Ec_X=\varnothing$. 
Hence, $\tau(b')\leq\theta$ for every $b'\in g^{-m}(a)$.
We deduce
from the definition of $\tau$ and $\theta$ that
$$\theta=\tau(a')\leq d^{-pm}\sum_{b'\in g^{-m}(a')} \tau(b')\leq
\theta.$$
It follows that $g^{-m}(a')\subset E$. Therefore, the analytic subset $\overline E$ 
of $Y_\epsilon$ satisfies $g^{-m}(\overline E)\subset\overline
E$. This contradicts the maximality of $\Ec_X$. 
\endproof

\noindent
{\bf End of the proof of Theorem \ref{th_conjecture_k}.}
Let $a$ be a point outside $\Ec_X$. Fix $\epsilon>0$ and a constant $\alpha>0$ small
enough. If $\nu$ is a limit value of $d^{-pmn}(g^{mn})^*(\delta_a)$, it is
enough to show that $\|\nu-s^{-1}\mu_X\|\leq 2\alpha+\epsilon$. Proposition 6.8 implies that
$\tau(a)=0$. So for $r$ large enough we have $\tau_r(a)\leq\alpha$. 
Consider all the negative orbits $\Oc_j$ of order $r_j\leq r$ 
$$\Oc_j=\{a^{(j)}_{-r_j},\ldots,a^{(j)}_{-1},a^{(j)}_{0}\}$$
with  $g^m(a^{(j)}_{-i-1})=a^{(j)}_{-i}$ and $a^{(j)}_{0}=a$ such that 
$a^{(j)}_{-r_j}\not\in Y_\epsilon$ and $a^{(j)}_{-i}\in Y_\epsilon$ for
$i\not=r_j$. Each orbit is repeated according to its multiplicity. 
Let $S_r$ denote the family of points $b\in
g^{-mr}(a)$ such that $g^{mi}(b)\in Y_\epsilon$ for $0\leq i\leq r$.
Then $g^{-mr}(a)\setminus S_r$ consists of the preimages of the points
$a^{(j)}_{-r_j}$. So, by definition of $\tau_r$, we have
$$d^{-pmr}\#S_r=\tau_r(a)\leq \alpha$$
and
$$d^{-pmr}\#(g^{-mr}(a)\setminus S_r)=d^{-pmr}\sum_j d^{pm(r-r_j)}=1-\tau_r(a)\geq 1-\alpha.$$
We have for $n\geq r$
\begin{eqnarray*}
d^{-pmn}(g^{mn})^*(\delta_a) & = & d^{-pmn}\sum_{b\in S_r} 
(g^{m(n-r)})^*(\delta_b) + d^{-pmn}\sum_j
(g^{m(n-r_j)})^*(\delta_{a^{(j)}_{-r_j}}).
\end{eqnarray*}
Since $d^{-pmn}(g^{mn})^*$ preserves the mass of any measure, the first term
in the last sum is of mass $d^{-pmr}\#S_r=\tau_r(a)\leq \alpha$ and the
second term is of mass $\geq 1-\alpha$. We apply Proposition \ref{prop_Briend_Duval} to
the Dirac masses at $a^{(j)}_{-r_j}$. We deduce that if $\nu$ is a limit value of
$d^{-pmn}(g^{mn})^*(\delta_a)$ then
$$\| \nu-s^{-1}\mu_X\|\leq 2\alpha +(1-\alpha)\epsilon\leq 2\alpha+\epsilon.$$
This completes the proof of the theorem.
\hfill $\square$
\par

\begin{corollary}
The cone of positive measures on $X$ which are totally invariant by
$g^m$, is of finite dimension.
In particular, the cone of positive measures on $\P^k$ which are totally invariant by
$f$, is of finite dimension.
\end{corollary}
\proof
Replacing $f$ by an iterate allows to assume that $g^m$ fixes all the
components of every analytic set which is totally invariant by $g^m$.
So, all these components are totally invariant.
Let $\nu$ be an extremal probability measure totally invariant by $g^m$. Let $X'$ be the smallest
analytic set totally invariant by $g^m$ such that $\nu(X')=1$. Since
$\nu$ is extremal, $X'$ is irreducible and $\nu(\Ec_{X'})=0$. It follows from Theorem
\ref{th_conjecture_k} and the invariance of $\nu$ that $\nu$ is
proportional to $\mu_{X'}$. By Corollary \ref{cor_excep_finite}, the
family of such measures is finite. 
\endproof

The following lemma will be useful in the proof of our main results
where $n_0$ is an index such that $\Ec^n_X=\Ec_X$ for $n\geq n_0$.

\begin{lemma} \label{lemma_pullback_variety}
There is a constant $\theta>0$ such that if $Z$ is an
  analytic subset of pure dimension $q\leq p-1$ of $X$ not contained
  in $\Ec_X$ then for every $n\geq 0$,
  $g^{-mn}(Z)$ contains an analytic set $Z_{-n}$ of pure dimension $q$
  of degree $\geq \theta d^{mn(p-q)}$. Moreover, if $n\geq n_0$ and if $x$ is a generic point in
  $Z_{-n}$, then $x\in\reg(X)$, $g^{m(n-n_0)}(x)\in\reg(X)$ and
  $g^{m(n-n_0)}$ defines a biholomorphism between a neighbourhood of
  $x$ and a neighbourhood of $g^{m(n-n_0)}(x)$ in $X$.
\end{lemma}
\proof
Let $P$ be a generic projective plane in $\P^k$ of dimension
$k-q$. Consider a point 
$a$ in $Z\cap P\setminus\Ec_X$. Since $\Ec_X=\Ec_X^{n_0}$, we have $g^{-ml}(a)\not\subset E_1$ for some $0\leq
l\leq n_0$. Then, by
 Lemma \ref{lemma_exceptional},
$\#g^{-mn}(a)$ contains
at least $\gamma d^{mp(n-n_0)}$ distinct points $x$ satisfying the last property in
the lemma. Let $Z_{-n}$ denote the union of the irreducible components of
$g^{-mn}(Z)$ which contain at least one such point $x$.
Then, $Z_{-n}$ satisfies the last property in the lemma.
We have $\#Z_{-n}\cap f^{-mn}(P) \geq \gamma d^{mp(n-n_0)}$. Since 
$\deg f^{-mn}(P)=d^{mnq}$, we obtain that $\deg Z_{-n}\geq  \theta
d^{m(p-q)n}$ for $\theta:=\gamma d^{-mpn_0}$.
\endproof


\section{Convergence towards the Green current} \label{section_cv}

In this section, we will prove the main results. 
Define the exceptional set $\Ec$ as the union of proper analytic
subsets $E$ of $\P^k$ which are totally invariant by $f$ and are {\it minimal}
in the following sense. The set $E$ does not contain non-empty proper
analytic sets which are totally invariant by $f$. 
Theorem \ref{th_exceptional} and 
Corollary \ref{cor_excep_finite} imply that $\Ec$ is a totally
invariant analytic set and it does not change if we replace $f$ by an
iterate of $f$, see also Remark \ref{rk_exception}.
We have the following result which implies Theorems \ref{main_theorem}
and \ref{main_theorem_bis}.

\begin{theorem} \label{th_final}
Let $f$, $T$, $\Ec$ be as above. Let $\Gc$ be a family of
modulo $T$ psh functions on $\P^k$ which is bounded in $L^1(\P^k)$. 
Assume that the restriction of $\Gc$ to
each component of $\Ec$ is a bounded family of modulo $T$ wpsh functions.
Then, $d^{-n}u\circ f^n$ converge to $0$ in $L^1(\P^k)$ uniformly on
$u\in\Gc$.  
\end{theorem}  

Let $m\geq 1$ be an integer such that $f^m$ fixes all the irreducible components
of all the totally invariant analytic sets. By Proposition \ref{prop_limit_lelong}, we can replace $f$
by $f^m$ and assume that $f$ fixes all these components.  
Let $X_p$ denote the
union of totally invariant sets of pure dimension $p$.
We will prove by induction on $p$ that 
 $d^{-n}u\circ f^n$ converge to $0$ in $L^1(X_p)$ uniformly on
$u\in\Gc$. We obtain the theorem for $p=k$ and $X_k=\P^k$.
Assume this convergence on $X_0$, $\ldots$, $X_{p-1}$ (the case $p=0$ is
clear). 
Define $X:=X_p$
and $\Ec_X$ as in Section \ref{section_exceptional}. From the
induction hypothesis, on each component $E$ of $\Ec_X$, 
 $d^{-n}u\circ f^n$ converge in $L^1$ to  0 uniformly 
on $u\in\Gc$. We deduce that $\Gc$ is bounded in $L^1(E)$. So,
if $Z$ is a
component of $X$, which is not minimal in the above sense,
by Lemma \ref{lemma_restriction_compact_T},  $\Gc$
is bounded in $L^1(Z)$. 
If $Z$ is a minimal component of $X$, then by hypothesis of the theorem, $\Gc$ is
bounded in $L^1(Z)$. So, we can apply Corollary \ref{cor_psh_convergence_T} to $\Gc$.

Let $\Gc'$ denote the set of all the modulo $T$ wpsh functions on $X$
which are limit values in $L^1(X)$ of a sequence $(d^{-n} u_n\circ f^n)$ with
$u_n\in\Gc$. For every $u\in \Gc'$, Corollary \ref{cor_psh_convergence_T} implies that $u\leq
0$. Since $\Ec_X\subset X$, by induction hypothesis we have convergence on $\Ec_X$ .
The last assertion of  Proposition \ref{prop_T_psh_compactness} 
imply that $u\geq 0$ on $\Ec_X$. Hence, $u=0$ on $\Ec_X$ 
for every $u\in\Gc'$.  It is clear that $\Gc'$ is compact. 
Fix a function $v_0\in\Gc'$. We have to show that $v_0=0$.

\begin{lemma} There are
  functions $v_n\in\Gc'$ such that $v_{n+1}=d^{-1}v_{n}\circ f$ almost
  everywhere for $n\in \Z$. 
\end{lemma}
\proof
Assume that $v_0$ is the limit of a sequence $(d^{-n_i}u_{n_i}\circ f^{n_i})$. Then,
for $n\geq 0$ the sequence $(d^{-n_i-n} u_{n_i}\circ f^{n_i+n})$ converges to
$d^{-n}v_0\circ f^n$. Lemma \ref{lemma_t_psh_interate} 
implies that $d^{-n}v_0\circ f^n$ is equal almost everywhere to an
element $v_n$ of $\Gc'$.
If $v_{-1}\in\Gc'$ is a limit value of $(d^{-n_i+1} u_{n_i}\circ f^{n_i-1})$
then $v_0=d^{-1}v_{-1}\circ f$ almost everywhere.
We construct the functions $v_{-n}$ in the same way by 
induction. If $v_{-n}$ is the limit of $(d^{-m_i}u'_{m_i}\circ f^{m_i})$ then 
we obtain $v_{-n-1}$ as a limit value of  $(d^{-m_i+1}u'_{m_i}\circ f^{m_i-1})$.
\endproof

\noindent
{\bf Proof of Theorem \ref{th_final}.} 
Let $\Gc''$ denote the set of all the modulo $T$ wpsh functions $w$ on $X$
which are limit values of the sequence $(v_{-n})_{n\geq 0}$.
Since $\Gc'$ is compact, we have $\Gc''\subset\Gc'$. 
We have to show that $v_0=0$. Assume this is not the case. 
Since $v_0=d^{-n}v_{-n}\circ f^n$ almost everywhere, 
by Lemma \ref{lemma_T_psh_lelong_zero},  there is a constant
$\alpha_0>0$ such that
$\max_X\nu_X(w,a)\geq\alpha_0$ for every $w\in\Gc''$. Fix a function $w_0\in\Gc''$.

\begin{lemma} There are
  functions $w_n\in\Gc''$ such that $w_{n+1}=d^{-1}w_{n}\circ f$ 
almost everywhere for $n\in \Z$. 
\end{lemma}
\proof
Assume that $w_0$ is the limit of $(v_{-n_i})$.
Let $w_1$ and $w_{-1}$ be modulo $T$ wpsh functions which are limit values of $(v_{-n_i+1})$ and
$(v_{-n_i-1})$ respectively. These functions belong to $\Gc''$. 
Then, $w_0=d^{-1}w_{-1}\circ f$ and $w_1=d^{-1}w_0\circ
f$ almost everywhere. We obtain the lemma by induction. If $w_n$ is the limit
value of $(v_{-m_i})$ then we obtain $w_{n-1}$ or $w_{n+1}$ as limit
values of $(v_{-m_i-1})$ or $(v_{-m_i+1})$ respectively.
\endproof

For $\alpha>0$ and $0\leq q\leq p-1$, denote by $\N_{\alpha,q}$
(resp. $\N_{\alpha,>q}$) the
family of indices $n\in\N$ such that $\{\nu_X(w_{-n},a)\geq \alpha\}$
is a non-empty analytic set of dimension $q$ (resp. $>q$). 
From the definition of $\Gc''$, we have
 $\cup_q \N_{\alpha_0,q}=\N$. Hence, there is a maximal
integer $q$ such that the upper density
$$\Theta^*(\N_{\alpha,q}):=\limsup_{n\rightarrow
  \infty} \frac{\#\N_{\alpha,q}\cap \{0,\ldots, n-1\}}{n}$$
is strictly positive for some constant $\alpha>0$. Fix a constant
$0<\beta\ll \alpha$ that we will choose later. The maximality of $q$
implies that $\Theta^*(\N_{\beta,>q})=0$. It follows that 
$$\delta:=\Theta^*(\N_{\alpha,q}\setminus \N_{\beta, >q})=\Theta^*(\N_{\alpha,q})>0.$$
Hence, for any integer $l\geq 1$, there is an integer $n_1\in
\N_{\alpha,q}\setminus \N_{\beta,>q}$ such that 
$$\#(\N_{\alpha,q}\setminus\N_{\beta,>q})\cap \{n_1,\ldots,n_1+l\}\geq
\frac{1}{2}\delta l.$$

Fix $l$ large enough and choose $\beta=\frac{1}{2} d^{-l-lk^2}\alpha$.
Replacing $w_0$ by $w_{-n_1}$
allows us to assume that $n_1=0$. This simplifies the notation.
We are looking for a contradiction using Lemma
\ref{lemma_demailly_meo} applied to $u:=w_0$. 
The hypothesis on the dimension of $\{\nu(w_0,a)>\beta\}$ is satisfied
since $0\in \N_{\alpha,q}\setminus \N_{\beta,>q}$. 
Let $n_0$  be given in Theorem \ref{th_exceptional} 
and Lemma \ref{lemma_pullback_variety}.
Choose integers $n_0<i_1<\cdots<i_s\leq l$, with $s\geq
\frac{1}{2} \delta l-n_0-1$,  in $\N_{\alpha,q}\setminus\N_{\beta,>q}$.
Let $Z_r''$ be an irreducible analytic set of
dimension $q$ such that $\nu_X(w_{-i_r},x)\geq \alpha$ on
$Z_r''$. 
We have seen that $w_{-i_r}=0$ on $\Ec_X$, hence $Z_r''\not\subset \Ec_X$.
By Lemma \ref{lemma_pullback_variety} (we assumed that $m=1$), there are analytic sets
$Z_r\subset g^{-i_r}(Z_r'')$
of pure dimension $q$ and of degree  $\geq \theta d^{i_r(p-q)}=:d_r$ such that
if $x$ is a generic point in $Z_r$ then $x\in\reg(X)$,
$x':=g^{i_r-n_0}(x)\in\reg(X)$ and $g^{i_r-n_0}$ defines a 
biholomorphism between neighbourhoods of $x$ and $x'$.
We now check the assumption of Lemma \ref{lemma_demailly_meo} that the
Lelong number of $w_0$  is $\geq
2\beta$ on $Z_r$.

Since
$w_0=d^{-i_r+n_0}w_{-i_r+n_0}\circ g^{i_r-n_0}$, we deduce from the
previous property of $g^{i_r-n_0}$ that 
$$\nu_X(w_0,x)=d^{-i_r+n_0}\nu_X(w_{-i_r+n_0},x').$$ 
Define $x'':=g^{i_r}(x)=g^{n_0}(x')$. This is a point in $Z_r''$. 
The local topological degree of $f^{n_0}$ at $x'$ is $\leq
d^{n_0k}$.  
Proposition
\ref{prop_lelong_compare} applied to $h:=f^{n_0}$ and the identity $w_{-i_r+n_0}=d^{-n_0}
w_{-i_r}\circ g^{n_0}$
imply that
$$\nu_X(w_{-i_r+n_0},x')\geq d^{-n_0-n_0k^2}\nu_X(w_{-i_r},x'')\geq
d^{-n_0-n_0k^2}\alpha.$$
It follows that $\nu_X(w_0,x)\geq d^{-i_r-n_0k^2}\alpha=:\nu_r\geq
2\beta$. Applying 
Lemma \ref{lemma_demailly_meo} yields 
$$\theta d^{-n_0k^2(p-q)}\alpha^{p-q} s\leq
2^{p-q+1}\deg(X)^{p-q}.$$ 
This is a contradiction if $l$ is large
enough, since $s\geq \frac{1}{2}\delta l -n_0-1$.
\hfill $\square$

\bigskip

\noindent
{\bf Proof of Theorem \ref{theorem_generic_map}.} It is enough to
prove that for $f$ generic in ${\cal H}_d$ we have $\Ec_{\P^k}=\varnothing$. By Lemma
\ref{lemma_topological} applied to $X=\Ec_{\P^k}$, it is enough to show that if $f$ is generic,
$\limsup d^{-(k-1)n}\#f^{-n}(x)=+\infty$ for every $x\in\P^k$. Here,
we count points without multiplicity. Fix an
$m\geq 1$ such that $d^m> 2^kk!$. We show for
$f$ generic that $\#f^{-m}(x)>d^{m(k-1)}$ for every $x\in\P^k$. This implies the result. Observe
that the family of such $f$ is a Zariski open set in ${\cal H}_d$. So,
it is enough to construct an $f$ satisfying this property. 

Choose a
rational map
$h:\P^1\rightarrow \P^1$ of degree $d$ such that $\#h^{-m}(x)\geq
\frac{1}{2}d^m$ for every $x\in \P^1$. To this end, it is enough to take a map
$h$ whose critical points are simple and have disjoint orbits.  
Now, construct the map $f$ using an idea of Ueda \cite{Ueda}. Let $\pi:\P^1\times
\cdots\times \P^1\rightarrow \P^k$ denote the canonical map which identifies
all points $(x_1,\ldots,x_k)$ with the points obtained by 
permutation of coordinates. If $\widehat f$ is the endomorphism of $ \P^1\times \cdots\times
\P^1$, $k$ times,  defined by
$\widehat f(x_1,\ldots,x_k):=(h(x_1),\ldots,h(x_k))$, then there is a
holomorphic map $f:\P^k\rightarrow\P^k$ of algebraic degree $d$ such that
$f\circ\pi=\pi\circ\widehat f$. We also have $f^m\circ \pi=\pi\circ
\widehat f^m$. 
Consider a point $x$ in $\P^k$ and a
point $\widehat x$ in $\pi^{-1}(x)$. We have $\pi^{-1}(f^{-m}(x))=
\widehat f^{-m}(\pi^{-1}(x))$. Hence,
$\#\pi^{-1}(f^{-m}(x))\geq \#\widehat f^{-m}(\widehat x)\geq 2^{-k}d^{mk}$.
Since $\pi$ has degree $k!$, we obtain $\#f^{-m}(x)\geq
\frac{1}{2^kk!} d^{mk}>d^{m(k-1)}$. This completes the proof.
\hfill$\square$

\begin{remark} \rm
Let $\Cc$ denote the compact convex set of totally invariant
$(1,1)$-currents of mass 1 on $\P^k$. Define an operator $\vee$ on $\Cc$.
If $S_1$, $S_2$ are elements of $\Cc$, write $S_i=T+\ddc u_i$ with
$u_i$ psh modulo $T$ on $\P^k$ such that $u_i\leq 0$ and $u_i=0$ on
$\supp(\mu)$, see Corollary \ref{cor_psh_convergence_T}. Define
$S_1\vee S_2:=T+\ddc \max(u_1,u_2)$. It is easy to check that $S_1\vee
S_2$ is an
element of $\Cc$. An element $S$ is said to be {\it minimal} if
$S=S_1\vee S_2$ implies $S_1=S_2=S$. It is clear that $T$ is not minimal
if $\Cc$ contains other currents.  A current of integration on a totally
invariant hypersurface is a minimal element.
\end{remark}

\begin{example} \label{example1}
 \rm
Let $[z_0:\cdots:z_k]$ denote the homogeneous coordinates of
$\P^k$ and $\pi:\C^{k+1}\setminus\{0\}\rightarrow\P^k$ 
the canonical projection. Consider the map $f[z_0:\cdots:z_k]:=[z_0^d:\cdots:z_k^d]$,
$d\geq 2$. The Green $(1,1)$-current $T$ of $f$ is
given by $\pi^*(T)=\ddc (\max_i\log|z_i|)$, see \cite{Sibony}, or equivalently
$T=\omega+\ddc v$ where 
$$v[z_0:\cdots:z_k]:=\max_{0\leq i\leq k}\log |z_i|
-\frac{1}{2}\log(|z_0|^2+\cdots+|z_k|^2).$$
The currents $T_i$ of integration on $(z_i=0)$ belong to $\Cc$ and 
 $T_j=T+\ddc u_j$ with $u_j:=\log|z_j|-\max_i\log|z_i|$. These
 currents are minimal. If
$\alpha_0$, $\ldots$, $\alpha_k$ are positive real numbers such that 
$\alpha:=1-\sum\alpha_i$ is positive, then $S:=\alpha T+\sum \alpha_i
T_i$ is an element of $\Cc$. We have $S=T+\ddc u$ with $u:=\sum
\alpha_i u_i$. The current $S$ is minimal if and only if 
$\alpha=0$. One can obtain other elements of $\Cc$ using the
operator $\vee$. One can also prove that $\Cc$ 
admits an infinite number of elements which are extremal in the cone
of positive closed $(1,1)$-currents. This implies that $\Cc$ has
infinite dimension. The elements of the set $\Ec$ in this case are
just the points $[0:\cdots:0:1:0:\cdots:0]$.
\end{example}

\section{Polynomial automorphisms}

The approach that we used above can be extended to other
situations. From now on we consider a polynomial automorphism
$f:\C^k\rightarrow \C^k$ of degree $\geq 2$ and its extension as a birational map on $\P^k$
that we also denote by $f$. Let $I_+$ and $I_-$ denote the
indeterminacy sets of $f$ and $f^{-1}$ respectively. These are the
analytic sets where $f$ and $f^{-1}$ are not defined; they are
contained in the hyperplane at infinity $L:=\P^k\setminus \C^k$. Assume that $f$ is {\it regular},
i.e. $I_+\cap I_-=\varnothing$.  We refer the reader to \cite{Sibony}
for the basic properties of regular automorphisms. There is an integer
$1\leq s\leq k-1$ such that $I_+$ and $I_-$ are irreducible analytic
sets of dimension $s-1$ and $k-s-1$ respectively.  We also have
$f(L\setminus I_+)=f(I_-)=I_-$ and $f^{-1}(L\setminus I_-)=f^{-1}(I_+)=I_+$. The maps
$f^n$ and $f^{-n}$ are also regular. 
The algebraic degrees $d_+$ and $d_-$ of $f$ and
$f^{-1}$ satisfy the relation $d_+^{k-s}=d_-^s$. 

The Green currents of bidegree $(1,1)$ associated to $f$ and $f^{-1}$
are denoted by $T_+$ and $T_-$. They are limits in the sense of
currents of $d_+^{-n} (f^n)^*(\omega)$ and $d_-^{-n}(f^n)_*(\omega)$
respectively. The current $T_+$ has locally continuous potentials
outside $I_+$, the current $T_-$ has locally continuous potentials
outside $I_-$. We also have $f^*(T_+)=d_+T_+$ and $f_*(T_-)=d_-T_-$.
We will consider the problem of convergence towards $T_+$,
the case of $T_-$ is obtained in the same way.

Let $g:X\rightarrow X$ denote the restriction of $f$ to $X:=I_-$. The
positive measure $\mu_X:=T_+^{k-s-1}\wedge [X]$ has positive mass. Since
$T_+$ is totally invariant, we have $g^*(\mu_X)=d_+^{k-s-1}\mu_X$. This
implies that $g$ has topological degree $d_+^{k-s-1}$. We construct as
above the families $X_0$, $\ldots$, $X_{k-s-1}$ of totally invariant
sets associated to $g$ with $X_{k-s-1}=I_-$. Let $\Ec_+$ denote the union of minimal
components in $\{X_0,\ldots, X_{k-s-1}\}$. 
We have the following result, see \cite{FornaessSibony1} for the case
of dimension 2.

\begin{theorem}
Let $S$ be a  positive closed $(1,1)$-current of
mass $1$ on $\P^k$. Assume that the local
potentials of $S$ are not identically equal to $-\infty$ on any
irreducible component of $\Ec_+$. Then, $d_+^{-n}(f^n)^*(S)$ converge to $T_+$. 
\end{theorem}

The proof follows the same lines as above. We will describe the
difference with the case of holomorphic endomorphisms and leave the
details to the reader. There is a neighbourhood $V$ of $I_+$ with
smooth boundary, which
can be chosen arbitrarily small,
such that $f(\P^k\setminus V)\Subset \P^k\setminus V$, see
\cite{Sibony}. If $S$ is as above, there is a modulo $T_+$
psh function $u$ such that $S=T_++\ddc u$. This function is defined
and is locally
bounded from above on $\P^k\setminus I_+$.
Denote by 
$\Gc$ the set of modulo $T_+$ psh
functions on $\P^k$ which are limit values of $d_+^{-n}u\circ f^n$.
Since the Lelong number of $u$ is $\leq 1$ at every point in
$\P^k\setminus I_+$ and since $f$ is an automorphism, Proposition
\ref{prop_lelong_compare} implies that the Lelong number of $d_+^{-n}u\circ f^n$ is $\leq
d_+^{-n}$ at every point in $\C^k$.

On the other hand, for $v\in\Gc$, we prove as in the previous sections that $v\leq 0$ and
$v=0$ on $X=I_-$. It follows that 
$v=0$ on $L\setminus I_+$ since we can write $v=d_+^{-1}v'\circ f$
with $v'\in\Gc$ and $f(L\setminus I_+)=I_-$.
The upper semi-continuity of the Lelong number implies that for every $\delta>0$, there is an $m$ such that the Lelong
number of $d_+^{-m}u\circ f^m$ is smaller than $\delta$ on $\P^k\setminus V$.  
We want to prove that $v=0$ on $\P^k\setminus \overline V$.

Assume  that $v=\lim
d_+^{-n_i} u\circ f^{n_i}$ and that
$v\leq -2\alpha$ with $\alpha>0$, on a ball
$B\subset\P^k\setminus\overline V$ of radius $r$.  Then as in Proposition
\ref{prop_limit_lelong}, we will have 
that $d_+^{-m}u\circ f^m\leq -d_+^{n_i-m}\alpha$ on a ball $B_i\subset
\P^k\setminus \overline V$ of radius
$\exp(-cr^{-2k}d_+^{n_i-m})$; this contradicts Proposition \ref{prop_psh_log} for
$\delta$ small and $n_i$ large. We can also obtain a uniform
convergence for 
regular automorphisms as in Theorem \ref{th_final}.


Tien-Cuong Dinh, UPMC Univ Paris 06, UMR 7586, Institut de
Math\'ematiques de Jussieu, F-75005 Paris, France. {\tt
  dinh@math.jussieu.fr}, {\tt http://www.math.jussieu.fr/$\sim$dinh} 

\

\noindent
Nessim Sibony,
Universit{\'e} Paris-Sud, Math{\'e}matique - B{\^a}timent 425, 91405
Orsay, France. {\tt nessim.sibony@math.u-psud.fr}

\end{document}